\documentclass[12pt]{article}

\usepackage{graphicx}
\usepackage{geometry}
\usepackage{hyperref}
\usepackage{amsmath,amssymb,esint,amsthm,bbm,textgreek,mathrsfs,mathtools}
\usepackage{setspace}
\usepackage{etoolbox}

\usepackage{stackengine}
\usepackage{scalerel}

\usepackage{cleveref}
\usepackage{algpseudocode}
\usepackage{algorithm}
\usepackage{caption}
\usepackage{subfigure}
\usepackage{graphicx}
\usepackage[normalem]{ulem}
\usepackage[table]{xcolor}
\usepackage{color, colortbl}
\usepackage{booktabs}
\usepackage{pdflscape}
\usepackage{tkz-fct}  
\usepackage{placeins}

\usepackage{pgfplots}
\pgfplotsset{compat=1.15}
\usepackage{mathrsfs}
\usetikzlibrary{arrows}

\usepackage{newtxtext}
\usepackage{natbib}
\hypersetup{
	colorlinks = true,
	urlcolor   = blue,
	citecolor  = black,
}

\newcommand{\RomanNumeralCaps}[1]
\linenumbers
\setcitestyle{square}

\crefname{equation}{}{}

\DeclareFontFamily{U}{mathx}{\hyphenchar\font45}
\DeclareFontShape{U}{mathx}{m}{n}{
      <5> <6> <7> <8> <9> <10>
      <10.95> <12> <14.4> <17.28> <20.74> <24.88>
      mathx10
      }{}
\DeclareSymbolFont{mathx}{U}{mathx}{m}{n}
\DeclareFontSubstitution{U}{mathx}{m}{n}
\DeclareMathAccent{\widecheck}{0}{mathx}{"71}
\makeatletter
\def\widebreve{\mathpalette\wide@breve}
\def\wide@breve#1#2{\sbox\z@{$#1#2$}%
     \mathop{\vbox{\m@th\ialign{##\crcr
\kern0.08em\brevefill#1{0.8\wd\z@}\crcr\noalign{\nointerlineskip}%
                    $\hss#1#2\hss$\crcr}}}\limits}
\def\brevefill#1#2{$\m@th\sbox\tw@{$#1($}%
  \hss\resizebox{#2}{\wd\tw@}{\rotatebox[origin=c]{90}{\upshape(}}\hss$}
\makeatletter
\newcommand\superwidetilde[1]{\stackon[-6pt]{$#1$}{\vstretch{1.3}{\hstretch{2.4}{\widetilde{\phantom{\;\;\;\;\;\;\;}}}}}}
\newcommand\ssuperwidetilde[1]{\stackon[-6pt]{$#1$}{\vstretch{1.3}{\hstretch{3.0}{\widetilde{\phantom{\;\;\;\;\;\;\;\;\;}}}}}}

\newcommand{\E}{\mathcal{E}}

\newcommand{\K}{\mathcal{K}}
\newcommand{\M}{\mathcal{M}}
\newcommand{\N}{\mathcal{N}}
\newcommand{\I}{\mathcal{I}}
\newcommand{\J}{\mathcal{J}}
\renewcommand{\L}{\mathcal{L}}
\renewcommand{\P}{\mathcal{P}}
\newcommand{\tP}{\widetilde{\mathcal{P}}}
\newcommand{\tE}{\widetilde{\mathcal{E}}}
\newcommand{\R}{\mathcal{R}}

\newcommand{\X}{\mathcal{X}}
\newcommand{\Y}{\mathcal{Y}}
\newcommand{\TM}{\mathcal{TM}}

\renewcommand{\a}{\alpha}

\newcommand{\w}{\omega}

\newcommand{\RR}{\mathbb{R}}

\newcommand{\NN}{\mathbb{N}}

\renewcommand{\vec}[1]{\boldsymbol{#1}}

\renewcommand{\u}{\vec{u}}

\newcommand{\x}{\mathbf x}

\newcommand{\grad}{\nabla}
\newcommand{\bfgrad}{\boldsymbol\nabla}
\newcommand{\bfDelta}{\boldsymbol\Delta}

\DeclareMathOperator*{\argmin}{arg\,min}
\DeclareMathOperator*{\codim}{codim}
\newcommand{\Id}{\operatorname{Id}}

\newcommand{\gradperp}{\bfgrad^{\perp}}

\newcommand{\dx}{\Delta x}
\newcommand{\dt}{\Delta t}
\newcommand{\intO}{\int_{\Omega}}

\newcommand{\nuc}{\widecheck{\nu}}

\newcommand{\nuN}{\nu_{N}}

\newcommand\STATE[1]{ s }
\newcommand\STATEin[1]{ s }

\newcommand{\tbw}{\breve{{\w}}}

\newcommand{\smax}{s_{\text{max}}}

\newcommand{\cphi}{\widecheck{\phi}}
\newcommand{\cmu}{\widecheck{\mu}}
\newcommand{\clambda}{\widecheck{\lambda}}

\newcommand{\tw}{\widetilde{\w}}

\newcommand{\tpsi}{\widetilde{\psi}}

\newcommand{\tbpsi}{\breve{\psi}}
\newcommand{\fw}{f_{\w}}

\newcommand{\gtw}{\bfgrad \tw}
\newcommand{\gtbw}{\bfgrad \breve{\w}}

\newcommand{\cS}{\mathcal{S}}

\def\Bmp#1{ \begin{minipage}{#1} }
\def\Emp{ \end{minipage} }
\def\Bmpc#1{ \begin{minipage}[c]{#1} }
\def\Bmpt#1{ \begin{minipage}[t]{#1} }
\def\Bmpb#1{ \begin{minipage}[b]{#1} }

\renewcommand{\arraystretch}{1.5}
\definecolor{Gray}{gray}{0.9}

\newcommand{\revtt}[1]{{\color{black}#1}}

\newcommand\strike[2]{\bgroup\markoverwith
{\textcolor{red}{\rule[.5ex]{2pt}{0.4pt}}}\ULon{#1} \revtt{#2}}

\makeatletter
\renewcommand*\env@matrix[1][\arraystretch]{%
  \edef\arraystretch{#1}%
  \hskip -\arraycolsep
  \let\@ifnextchar\new@ifnextchar
  \array{*\c@MaxMatrixCols c}}
\makeatother

\renewcommand{\vec}[1]{\boldsymbol{#1}}

\newcommand\p[2]{\frac{\partial #1}{\partial #2}}

\newcommand{\ub}{\bar{u}_{b}}

\newcommand{\Nv}{\N_{\varphi}}

\newtheorem{problem}{Problem}

\newcommand{\cvarphi}{\widecheck{\varphi}}

\newcommand{\showprobC}[1]{Problem~\theproblem.#1}

\usepackage[T1]{fontenc}
\usepackage[utf8]{inputenc}
\usepackage{authblk}

\begin{document}
  \title{{Adjoint-Based Enforcement of {State} Constraints in PDE
    Optimization Problems}}

\author[1,2]{Pritpal Matharu\thanks{Email address for correspondence: \texttt{pritpal@kth.se}}}
\author[2]{Bartosz Protas}
\affil[1]{Department of Mathematics, KTH Royal Institute of Technology, Stockholm, Sweden}
\affil[2]{Department of Mathematics and Statistics, McMaster University, Hamilton, Ontario, Canada}


\date{July 25, 2024}
\maketitle

\begin{abstract}
  {This study demonstrates how the adjoint-based framework
    traditionally used to compute gradients in PDE optimization
    problems can be extended to handle general constraints on the
    state variables. This is accomplished by constructing a
    projection of the gradient of the objective functional onto a
    subspace tangent to the manifold defined by the constraint. This
    projection is realized by solving an adjoint problem defined in
    terms of the same adjoint operator as used in the system employed
    to determine the gradient, but with a different forcing. We focus
    on the ``optimize-then-discretize'' paradigm in the
    infinite-dimensional setting where the required regularity of both
    the gradient and of the projection is ensured. The proposed
    approach is illustrated with two examples: a simple test problem
    describing optimization of heat transfer in one direction and a
    more involved problem where an optimal closure is found for a
    turbulent flow described by the Navier-Stokes system in two
    dimensions, both considered subject to different state
    constraints.  The accuracy of the gradients and projections
    computed by solving suitable adjoint systems is carefully verified
    and the presented computational results show that the solutions of
    the optimization problems obtained with the proposed approach
    satisfy the state constraints with a good accuracy, although not
    exactly.}
\end{abstract}

\begin{flushleft}
\end{flushleft}



\section{Introduction \label{sec:Introduction}}

Numerous problems in modern science and engineering are cast in terms
of optimization of systems described by partial differential equations
(PDEs), often subject to complex constraints imposed on the control
(decision) variables. A standard framework for handling such
constrained PDE optimization problems is based on Lagrange multipliers
\cite{l69,Troltzsch2010} where the objective functional is considered
a function of both the state and the control variable, the two related
by suitable constraints. However, this formalism {can lead} to
challenging computational problems since the system of the first-order
optimality conditions has the form of a saddle-type problem which is
often difficult to solve numerically. As an alternative, one can
consider formulations based on the reduced objective functional
depending on the control variable only {(i.e., where the state
  variable is {\em not} a separate independent variable, but is
  considered a function of the control variable),} such that the
optimization problem can be solved using some discrete form of the
gradient flow in which the constraints are built into the definition
of the space in which the gradient flow is constructed. An advantage
of this approach is that a saddle-type (min-max) problem is avoided
and one can take advantage of computational methods of unconstrained
optimization. When dealing with infinite-dimensional PDE optimization
problems a key element of this approach is determination of the
gradient of the reduced objective functional with respect to the
control variable which can be conveniently computed by solving a
suitably-defined {\em adjoint system} \cite{control:lions1,g03}. These
methods are already quite well developed and there exists a large body
of literature on this topic{, with applications in various areas
  \cite{bp11a,Vergnault2014,Costanzo2022,Sirignano2023}}. However,
methods based on discrete gradient flows can be challenging to apply
in the presence of complicated constraints, especially when they
involve both state and control variables as is the case in some
applications {(one research area where such problems are common
  is topology optimization \cite{bendsoe2003topology}).} In the
present study we demonstrate how such constraints can be approximately
enforced using the solution of a suitable adjoint problem, similar but
distinct from the adjoint system used to determine the gradient of the
reduced objective functional.  We focus on the
``optimize-then-discretize'' approach where all elements of the
algorithm are derived in the continuous setting and only then
discretized for the purpose of the numerical solution \cite{g03}.  To
the best of our knowledge, this is the first such generalization of
adjoint-based techniques in the solution of PDE optimization problems.

To illustrate the ideas described above, we consider the state
$u \in \X$ and the control variable ${\varphi} \in \Y$ with $\X$ and
$\Y$ denoting the appropriate functional spaces, to be specified
below. The objective functional then is
$j \; : \; \X \times \Y \rightarrow \RR$ and the PDE constraint is
given in terms of some function
${S} \; : \; \X \times \Y \rightarrow \X^*$ representing the weak form
of the PDE with $\X^*$ the dual space of $\X$. In addition, we will
also assume that the state and the control variables are subject to
$m \in \NN^+$ scalar constraints given by {the function
  $c \; : \; \X \times \Y \rightarrow \RR^m$ assumed sufficiently
  smooth.} Our constrained PDE optimization problem can then be stated
as
\begin{equation*}
  \textrm{(A)} \hfill \qquad\qquad\qquad
  \begin{aligned}
    & \min_{(u,{\varphi}) \in \X \times \Y} j(u,{\varphi}) \\
    & \textrm{subject to}: \quad
    \left\{
    \begin{aligned}
      {S}(u,{\varphi}) &= 0, \\ c(u,{\varphi}) & = 0
    \end{aligned}\right. . 
  \end{aligned}
 \end{equation*}
 Introducing the Lagrange multipliers $\lambda \in \X$,
 $\mu \in \RR^m$ and the augmented objective functional, i.e., the
 Lagrangian,
 $\L(u,{\varphi},\lambda,\mu) := j(u,{\varphi}) + \left\langle \lambda,
   {S}(u,{\varphi}) \right\rangle_{\X \times \X^*} + \langle \mu, c(u,{\varphi})
 \rangle_{\RR^m}$, where
 $ \langle \cdot,\cdot \rangle_{\X \times \X^*}$ is the duality
 pairing between the spaces $\X$ and $\X^*$ \cite{l69},
 $\langle \cdot,\cdot\rangle_{\RR^m}$ the inner product in $\RR^m$,
 and ``:='' means ``equal to by definition'', the constrained problem
 (A) can be recast in the corresponding unconstrained form as
 \cite{l69}
\begin{equation*}
  \textrm{(B)} \hfill \qquad\qquad\qquad
 \min_{(u,{\varphi}) \in \X \times \Y} \quad  \max_{(\lambda,\mu) \in \X^* \times
   \RR^m} \L(u,{\varphi},\lambda,\mu).
\end{equation*}
Critical points of this problem are characterized by the appropriate
first-order optimality conditions. Besides handling the saddle-type
(min-max) nature of these critical points, an extra complication is
the need to also approximate the optimal Lagrange multipliers
$\clambda$ and $\cmu$ in addition to the optimal state and control
variables $\widecheck{u}$ and $\cphi$ {(hereafter, the symbol
  $\widecheck{\cdot}$ will denote the optimal value of a variable).}

Defining the {\em reduced} objective functional
$\J({\varphi}) := j(u({\varphi}),{\varphi})$, which assumes that for each ${\varphi}$ the
map $u = u({\varphi})$ is well defined in terms of the solution of the PDE
problem ${S}(u({\varphi}),{\varphi}) = 0$, Problem (A) can be recast such that
minimization is performed with respect to ${\varphi}$ only, i.e.,
\begin{equation*}
  \textrm{(C)} \hfill \qquad\qquad\qquad
  \begin{aligned}
    & \min_{{\varphi} \in  \Y} \J({\varphi}) \\
    & \textrm{subject to}: \quad  c(u({\varphi}),{\varphi})  = 0
  \end{aligned},
\end{equation*}
where the constraint can be interpreted as defining a codimension-$m$
manifold ${\M} := \left\{ {\varphi} \in \Y: \quad  c(u({\varphi}),{\varphi})  = 0
\right\}$. Thus, since the optimizer ${\widecheck{\varphi}}$ must be sought on this
constraint manifold, ${\widecheck{\varphi}} \in {\M}$, problem (C) defines a {\em
  Riemannian} optimization problem \cite{ams08}.

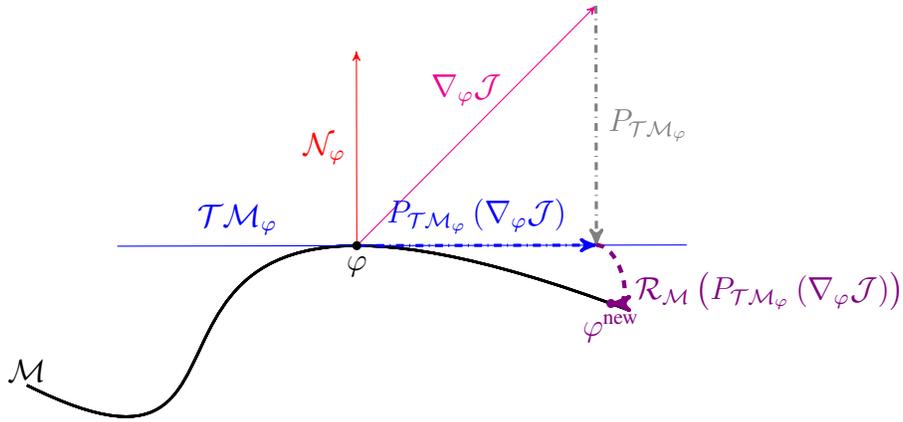
\begin{figure}
\centering
 \begin{tikzpicture}
\draw[black,line width=1pt]  (0,8) .. controls ++(4,-2) and ++(-8,3) .. (8,9)
      node[sloped,inner sep=0cm,above,pos=.792,
      anchor=south west,
      minimum height=(10.5)*0.3cm,minimum width=(10.5)*.3cm](M){};
\draw[black,line width=1pt]  (0,8) .. controls ++(4,-2) and ++(-8,3) .. (8,9)
      node[sloped,inner sep=0cm,above,pos=.792,
      anchor=south west,
      minimum height=(8.5)*0.3cm,minimum width=(14.5)*.3cm](N){};
\draw[black,line width=1pt]  (0,8) .. controls ++(4,-2) and ++(-8,3) .. (8,9)
      node[sloped,inner sep=0cm,above,pos=.792,
      anchor=north east,
      minimum height=(14.5)*0.3cm,minimum width=(10.5)*.3cm](P){};
\draw[black,line width=1pt]  (0,8) .. controls ++(4,-2) and ++(-8,3) .. (8,9)
      node[sloped,inner sep=0cm,above,pos=1-0.01,
      anchor=north east,
      minimum height=(14.5)*0.3cm,minimum width=(10.5)*.3cm](R){};

\path (M.south west)%
           edge[-stealth',magenta] node[left,pos=.65] {$\grad_{{\varphi}}\J$} (M.north east);
\path (M.south west)%
           edge[-stealth',red] node[left] {$\Nv$} (N.north west);
\draw[blue]  (M.south west)  --%
      node[above] {} (N.south east);
\draw[blue]  (P.north east)  --%
      node[above] {${\TM_{\varphi}}$} (P.north west);
\draw[-stealth',dashdotted, very thick, gray]  (M.north east)  --%
      node[right] {${P_{\TM_{\varphi}}}$} (M.south east);
\draw[-stealth',dashdotted,blue, very thick]  (M.south west)  --%
      node[above] {${P_{\TM_{\varphi}}}\left(\grad_{{\varphi}}\J\right)$} (M.south east);
\draw [-stealth', dashed, ultra thick, violet] (M.south east) to[out=1,in=0] node[right,pos=0.7] {$\R_{\M}\left({P_{\TM_{\varphi}}}\left(\grad_{{\varphi}}\J\right)\right)$} (R.north east);

\node[] at (0,8.2)   {$\M$};
\filldraw[black] (M.south west) circle (1.5pt) node[below] {$\varphi$};
\filldraw[violet] (R.north east) circle (1.5pt) node[below] {$\varphi^{\text{new}}$};
  \end{tikzpicture}    
\vspace{-1.2cm}
\caption{{Schematic representation of the projection of the gradient
    ${\grad_{\varphi}\J}$ onto the subspace $\TM_{{\varphi}}$ tangent
    to a general manifold $\M$ at ${\varphi}$. The projection operator
    ${P_{\TM_{\varphi}}}$ uses the element $\Nv$ normal to the tangent
    subspace $\TM_{{\varphi}}$, cf.~\cref{eq:PTM}. The retraction
    operator {$\R_{\M}$} {then maps elements of} the tangent
    subspace $\TM_{{\varphi}}$ back to the manifold $\M$,
    cf.~\cref{eq:retract}, to a new {element}
    ${\varphi}^{\text{new}}$}. {The objects shown in the figure
    are to be interpreted as elements of a function space, typically a
    Sobolev space \cite{af05}, and can be functions of time as in
    Problems \ref{prob:Heat1} and \ref{prob:Heat2}, or state as in
    Problem \ref{prob:nu}.}}
\label{fig:M}
\end{figure}

 A local {minimizer} of Problem (C) can then be approximated with a
 discrete projected gradient flow as ${\widecheck{\varphi}} = \lim_{n \rightarrow
   \infty} {\varphi}^{(n)}$, where
 \begin{subequations}
\label{eq:desc}
\begin{align}
{{\varphi}}^{(n+1)} &= {{\varphi}}^{(n)} - {\tau_{n}} \, {P_{\TM_{{\varphi}^{(n)}}}\left(\grad_{\varphi}\J({\varphi}^{(n)})\right)}, \quad\quad n=0, 1, \dots, \label{eq:desc1} \\
{{\varphi}}^{(0)} &= {{\varphi}}_0,
\end{align}
\end{subequations}
in which
${P}_{\TM_{{\varphi}}} \: : \; \Y \rightarrow \TM_{{\varphi}}$ is the
operator representing the {orthogonal} projection onto the
subspace $\TM_{{\varphi}} \in \Y$ tangent to the constraint manifold
$\M$ at ${\varphi} \in \M$, $\tau_n$ is the step size along the
projected gradient{, $\grad_{{\varphi}}\J$ is the gradient of the
  reduced functional {$\J(\varphi)$} with respect to
  ${\varphi}$,} and ${\varphi}_0$ the initial guess. As illustrated in
Figure \ref{fig:M}, in iterations \eqref{eq:desc} the constraint is
satisfied approximately only, with an error $\mathcal{O}(\tau_n^2)$ at
each iteration. This error is eliminated if one can introduce the
retraction operator ${\R_{\M}} \; : \: \TM \rightarrow \M$, where
$\TM$ is the tangent bundle, such that the right-hand side (RHS) in
\eqref{eq:desc1} is replaced with
${\R_{\M}\left({\varphi}^{(n)} - \tau_{n} \,
    P_{\TM_{{\varphi}^{(n)}}}\left(\grad_{\varphi}\J({\varphi}^{(n)})\right)
  \right)}$ yielding the Riemannian gradient approach
\cite{ams08}. However, for general manifolds $\M$ the retraction
operator can be very difficult to implement and as a result one often
needs to resort to using \eqref{eq:desc}.  While the gradient
${\grad_{{\varphi}}\J({\varphi})}$ of the (reduced) objective
functional can be conveniently computed by solving an adjoint problem
\cite{g03}, the main contribution of the present study is to
demonstrate that an analogous computation can in fact be also used to
determine the projection operator ${P}_{\TM_{{\varphi}}}$. This
{problem} is nuanced by the fact that both the gradient
${\grad_{{\varphi}}\J({\varphi})}$ and the projector
${P}_{\TM_{{\varphi}}}$ need to be determined with respect to a
particular topology, typically induced by the norm in the space
$\Y$. This contribution thus expands the scope of applications of
adjoint analysis in PDE constrained optimization.

We introduce and illustrate the new approach by way of two examples:
the first one is rather academic and concerns a simple control problem
involving a one-dimensional (1D) heat equation; the second is more
involved and represents an extension of a recently formulated problem
concerning finding optimal turbulence closures in two-dimensional (2D)
Navier-Stokes flows \cite{Matharu2022a}. While the first example is
intended to serve as {a simple illustration only of the proposed
  approach,} the second one showcases an application to a nontrivial
problem with a nonstandard structure.  The structure of the paper is
as follows: in the next section we introduce the two PDE optimization
problems whereas the {proposed} approach is presented in
\Cref{sec:solution} with a particular emphasis on how the projection
onto the subspace tangent to the constraint manifold can be
constructed based on the solution of an adjoint system; computational
results are then presented in \Cref{sec:results} with a summary and
conclusions deferred to \Cref{sec:final}.

\section{Test Problems \label{sec:prob}}

In this section we introduce two constrained PDE optimization problems
that will serve as our test cases. While the first problem is
classical and is used primarily to illustrate our approach, the second
one is more challenging due to its non-standard structure. With a
slight abuse of notation motivated by the desire to highlight the
analogies between them, in the two test problems we reuse {the
  symbols} denoting key elements of the formulation of the
optimization problems, with the exception of {$\phi$ and
  $\varphi$ which represent the control variable} in \Cref{sec:heat}
and \Cref{sec:closure}, respectively.

\subsection{Control of Heat Conduction in 1D}
\label{sec:heat}

We consider heat conduction in a finite rod
$\Omega := [a,b] \subset \RR$, where $-\infty < a < b < \infty$, in
which the time-dependent heat flux $\phi(t)$, $t \in [0,T]$ with some
$T > 0$, is applied at the left boundary $x = a$ such that the
resulting temperature at the right boundary $x = b$ should match some
prescribed history $\ub(t)$, $t \in [0,T]$, and {at the same
  time} the average ``energy'' of the system should be given by
$E_0 > 0$. The problem is thus described by the system
\begin{subequations} \label{eq:Heat}
\begin{alignat}{3} 
&\p{u}{t} - \Delta u& &= 0,& \qquad (t, x) &\in (0, T] \times \Omega, \\
&\p{u}{x}\Big|_{x=a}& &= \phi(t),& \quad t &\in (0, T], \label{eq:HeatBC} \\
&\p{u}{x}\Big|_{x=b}& &= 0,& \quad t &\in (0, T], \\
&u(t=0)& &= u_{0},& \quad x &\in \Omega, 
\end{alignat} 
\end{subequations}
where {$u_0 \; : \; \Omega \rightarrow \RR$ is the initial
  condition and} the solution $u = u(t,x;\phi)$ is assumed to depend
on the flux $\phi$ as the control variable. The energy {$E$} is
defined as
\begin{equation}
  \label{eq:E1}
E(t; \phi) := \frac{1}{2} \int_{\Omega} u(t, x; \phi)^2 \, dx,
\end{equation}
whereas $[f]_T := (1/T) \int_0^T f(t)\, dt$ will denote the time
average of some function $f \; : \; [0,T] \rightarrow \RR$. Next, we
introduce the Sobolev space $H^p(\Lambda)$ with $p \in \mathbb{N}$ of
functions defined on {some} domain $\Lambda$ with $p$
square-integrable distributional derivatives {and} endowed with
the inner product
\begin{equation}
\forall \  z_1,z_2 \in H^p(\Lambda) \qquad
\left\langle z_1, z_2 \right\rangle_{H^p(\Lambda)} 
{:=} \sum_{q=0}^{p} \ell_q^{2q}\left\langle \frac{d^q z_1}{dx^q}, \frac{d^q z_2}{dx^q} \right\rangle_{L^2(\Lambda)}, 
\label{eq:ipHp}
\end{equation}
where
$\langle z_1, z_2 \rangle_{L^2(\Lambda)} := \int_{\Lambda} z_1 z_2 \,
dx$, whereas $\ell_{q} \in \mathbb{R}^+$, $q = 0,\dots,p$, are
``length-scale'' parameters {(in the multidimensional case, i.e.,
  when $\Lambda \subseteq \RR^n$, $n>1$, the definition
  \eqref{eq:ipHp} admits a natural generalization \cite{af05})}. While
for different values of {$0 < \ell_q < \infty$} the inner
products in \eqref{eq:ipHp} are equivalent (in the precise sense of
norm equivalence), in \Cref{sec:results} we will show that the
flexibility represented by these parameters significantly improves the
convergence of the our iterative minimization algorithm
\eqref{eq:desc}.

{It is instructive to consider two formulations: one in which the
  flux is simply assumed to be a square-integrable function of time
  together with its first derivative, i.e.,
  $\phi \in \cS := H^1(0,T)$, and another one where it is additionally
  assumed to belong to a quadratic manifold,
  $\phi \in \M \subset \cS$, which is defined as}
\begin{equation} \label{eq:M}
  {\M} := \left\{ \phi \in {\cS}: \quad
    \left[E(\cdot; \phi)\right]_T = \frac{1}{2T} \int_0^T \int_{\Omega} u(t, x; \phi)^2 \, dx dt = E_0 \right\}.
\end{equation}
Our {(reduced)} objective functional
${\J} \; : \; {\M} \rightarrow \RR$ is a measure of the error between
the actual temperature at the right boundary $u(t,b;\phi)$ and the
prescribed profile $\ub(t)$
\begin{align} \label{eq:J}
{\J}(\phi) = \frac{1}{2} \int_{0}^T [u(\phi)|_{b} - \ub]^2 \, dt.
\end{align}
We note that the quadratic constraint defining the manifold ${\M}$ in
\eqref{eq:M} simplifies when the governing system \eqref{eq:Heat} is
subject to the homogeneous initial condition $u_0(x) = 0$,
$\forall \ x \in \Omega$ {($u_0 \equiv 0$)}. In such case the
quadratic constraint becomes homogeneous, i.e.,
$\left[E(\cdot; {\beta}\phi)\right]_T = {\beta}^2 \left[E(\cdot;
  \phi)\right]_T$ for some {$\beta \in \RR$,} which reduces
enforcement of this constraint to a simple rescaling. Otherwise,
enforcement of the constraint is more involved and necessitates the
use of a retraction operator.  {For comparison,} we thus consider
two versions of {each of the optimization problems} depending on
whether or not the initial condition $u_0$ is homogeneous, i.e.,
\begin{problem}\label{prob:Heat1}
Given system \cref{eq:Heat} with {$u_0 \equiv 0$}  and objective functional \cref{eq:J}, find
\begin{align} 
\cphi &= \underset{\phi \in {\cS}} {\argmin} \, {\J}(\phi), \label{prob:Heat1S} \tag{\showprobC{A}} \\
\cphi &= \underset{\phi \in {\M}} {\argmin} \, {\J}(\phi). \label{prob:Heat1M} \tag{\showprobC{B}}
\end{align}
\end{problem}

\begin{problem}\label{prob:Heat2}
Given system \cref{eq:Heat} with {$u_0 \not\equiv 0$}  and objective functional \cref{eq:J}, find
\begin{align}
\cphi &= \underset{\phi \in {\cS}} {\argmin} \, {\J}(\phi), \label{prob:Heat2S} \tag{\showprobC{A}} \\
\cphi &= \underset{\phi \in {\M}} {\argmin} \, {\J}(\phi). \label{prob:Heat2M} \tag{\showprobC{B}}
\end{align}
\end{problem} 
We add that since in \ref{prob:Heat1S} and \ref{prob:Heat2S}
minimization is performed over the entire linear subspace $\cS$, we
will consider these problems as ``unconstrained''. In contrast,
\ref{prob:Heat1M} and \ref{prob:Heat2M}, where optimization is carried
out over the manifold $\M \subset \cS$, will be regarded as
``constrained''.

\subsection{Optimal Eddy Viscosity in Closure Models {for} 2D Turbulent Flows}
\label{sec:closure}

Turbulent flows governed by the incompressible Navier-Stokes system
remain very challenging to compute when the Reynolds number is high
due to difficulties resolving motions occurring at small spatial and
temporal scales. A solution often adopted in practice relies on
solving a filtered version of the governing equations such that the
number of resolved degrees of freedom is reduced, an approach referred
to as the Large-Eddy Simulation (LES)
\citep{Lesieur1993book,pope2000turbulent,davidson2015turbulence}.
However, the LES system is not closed as it involves terms
{explicitly} depending on the {unresolved} degrees of
freedom and closing this system in terms of the resolved degrees of
freedom leads to the celebrated turbulence closure problem. While a
significant body of literature has been devoted to deriving turbulence
models adapted to different flows, most of these approaches have been
empirical in nature. Recent advances in machine learning have enabled
the development of data-based techniques for deducing closure models
\citep{kutz_2017,GamaharaHattori2017,jimenez_2018,duraisamy2018turbulence,Duraisamy2021,San2021,Waschkowski2022}.
As an alternative to these approaches, a data-based technique relying
on calculus of variations and methods of PDE optimization was recently
proposed in {\cite{Matharu2020,Matharu2022a}}. It allows one to find
an optimal, in a mathematically precise sense, functional form of the
eddy viscosity appearing in a commonly used family of turbulence
closure models. The problem considered here arises as an extension of
that approach.

We consider the flow of viscous incompressible fluids on a
two-dimensional (2D) periodic domain $\Omega = [0,2\pi ]^2$ and over a
time interval $[0,T]$ for some $T > 0$. Assuming uniform fluid density
$\rho = 1$, the motion is governed by the Navier-Stokes system, which
written in vorticity-streamfunction form, is
\begin{subequations} \label{eq:2DNS}
	\begin{alignat}{2} 
	\partial_t {w} + \gradperp \psi \cdot \bfgrad {w} &= \nuN \bfDelta {w}  - \a {w} + \fw & \quad &\text{in} \quad {(0, T] \times \Omega},  \label{eq:vort} \\
	\bfDelta \psi &= -{w} & \quad &\text{in} \quad {(0, T] \times \Omega}, \label{eq:stream_vort} \\
	{w}(t=0) &= {w_0} & \quad &\text{in} \quad \Omega, \label{eq:IC}
	\end{alignat} 
\end{subequations}
where {$w := -\gradperp \cdot \u$}, with
$\gradperp := [\partial_{x_2}, -\partial_{x_1}]^{T}$ and $\u$ the
velocity field, is the vorticity component perpendicular to the plane
of motion, $\psi$ is the streamfunction, $\nuN$ is the coefficient of
the kinematic viscosity (for simplicity {of notation}, we reserve
the symbol $\nu$ for the eddy viscosity), and ${w_0}$ is the initial
condition. In \eqref{eq:vort} the term proportional to $\alpha > 0$
represents large-scale dissipation ({Ekman} friction), whereas $\fw$
is a time-independent band-limited forcing acting on low
wavenumbers. The parameters of these two terms are chosen such that
the flow is in a statistical equilibrium with a well-developed
enstrophy cascade and a rudimentary inverse energy cascade
\citep{Lesieur1993book}.

Using $\widetilde{(\cdot)}$ to denote a suitable low-pass filter, we
can write the LES version of \cref{eq:2DNS} as (the reader is referred
to \cite{Matharu2022a} for derivation details)
\begin{subequations} \label{eq:LES}
	\begin{alignat}{2} 
	\partial_t \tw + {\superwidetilde{\gradperp \tpsi \cdot \bfgrad \tw}} &= \bfgrad \cdot {\ssuperwidetilde{\left( \left[\nuN + \nu(s) \right] \bfgrad \tw \right)}}  - \a \tw + \fw & \quad &\text{in} \quad {(0, T] \times \Omega},  \label{eq:LES_eqn} \\
	\bfDelta \tpsi &= -\tw & \quad &\text{in} \quad {(0, T] \times \Omega}, \label{eq:LES_stream_vort} \\
	\tw(t=0) &= \tw_0 {:=} \widetilde{w}_0& \quad &\text{in} \quad \Omega, \label{eq:LES_IC}
	\end{alignat} 
\end{subequations}
where $\tw$ is the LES {vorticity,} and the initial condition
{$\tw_0$} is given as the filtered initial condition
\eqref{eq:IC} and the forcing term is unaffected by the filter as it
acts on the low wavenumbers only. The LES equation \eqref{eq:LES_eqn}
features a Smagorinsky-type closure model with a state-dependent eddy
viscosity expressed as
\begin{equation}
  \nu(s) = \left[ {\eta^3} \, \sqrt{s} + \nu_0 \right] \varphi\left( \frac{s}{\smax} \right) \qquad \text{with} \quad s :=
|\gtw|^2 \in \I := [0,\smax],
  \label{eq:nu}
\end{equation}
where $\nu_0 > 0$, {$\eta = 2\pi/k_c$} is the width of the LES filter
with $k_c$ the largest resolved wavenumber, $\smax > 0$ is a
sufficiently large number to be specified later and
$\varphi \; : \; [0,1] \rightarrow \RR$ a non-dimensional
function. The form of equation \eqref{eq:LES_eqn} suggests that
$\nu = \nu(s)$, and hence also $\varphi = \varphi( s / \smax )$, must
be at least piecewise $C^1$ functions on $\I$ and $[0,1]$,
respectively.  However, as will become evident in {Section}
\ref{sec:nu}, our solution approach imposes some additional regularity
requirements, namely, $\nu = \nu(s)$ needs to be piecewise
{$C^3$} on $\I$ with the first and third derivatives vanishing at
$s = 0,\smax$. Since gradient-based solution approaches to
PDE-constrained optimization problems are preferably formulated in
Hilbert spaces \citep{pbh04}, we shall look for an optimal function
$\varphi$ parametrizing the eddy viscosity as an element of the
following linear space which is a subspace of the Sobolev space
$H^2(\I)$, cf~\eqref{eq:ipHp},
\begin{equation}
\cS := {\left\{ \varphi \in {C^3([0,1])} \: : \: \frac{d}{d\xi} \varphi(\xi) =  \frac{d^3}{d\xi^3} \varphi(\xi) = 0 \ \text{at} \ \xi = 0,1 \right\}}.
\label{eq:S}
\end{equation}
We will seek an optimal form of the function $\varphi$ such that
solutions of the LES system \eqref{eq:LES} with eddy viscosity
\eqref{eq:nu} best match the corresponding solutions of the original
Navier-Stokes system \eqref{eq:2DNS} in a sense to be specified below.
We add that when $\nu_0 = 0$ and $\varphi(\xi) = 1$, $\xi \in [0,1]$,
then the eddy viscosity \eqref{eq:nu} reduces to the Leith model which
is a counterpart of the Smagorinsky model in 2D flows
\citep{Leith1968, Leith1971, Leith1996}.

While in \cite{Matharu2022a} we considered a pointwise match, both in
space and in time, between the LES and the Navier-Stokes flows, here
we determine an optimal eddy viscosity $\nuc(s)$ that will allow the
LES flow to match the original Navier-Stokes flow in a certain average
sense. {{More specifically, we will formulate the problem
    in terms of time-averages involving important integral} quantities
  characterizing 2D flows, {namely, the enstrophy and
    palinstrophy which are} defined as follows}
\begin{subequations}
\label{eq:Ep}
\begin{align}
\E(t) & := \frac{1}{2} \intO \widetilde{w}(t, \x)^2 \, d\x, \label{eq:E} \\
\P(t) & := \frac{1}{2} \intO \left|\bfgrad \widetilde{w}(t, \x) \right|^2 \, d\x, \label{eq:P}
\end{align}
\end{subequations}
where $\widetilde{w}$ is the filtered solution of the Navier-Stokes
system \eqref{eq:2DNS}, whereas $\tE(t;\varphi)$ and $\tP(t;\varphi)$
will denote the corresponding quantities defined in terms of the
solution $\tw$ of the LES system \eqref{eq:LES} and the notation
emphasizes their dependence on the function $\varphi$ parameterizing
the eddy viscosity, cf.~ansatz \eqref{eq:nu}.

The optimal eddy viscosity $\nuc(s)$ will be chosen {so as to}
yield a good match between the palinstrophy in the LES flow (with the
given eddy viscosity) and in the original Navier-Stokes flow, so that
we shall consider the following error functional 
 \begin{equation}
\label{eq:J2}
\J({{\varphi}}) := \frac{1}{2\, D} \int_0^T \left[ \P(t) - \tP(t; {\varphi}) \right]^2 \, dt, 
\end{equation}
where {$D>0$} serves as a normalization factor. In addition, we
will also require the LES flow with the optimal eddy viscosity to have
a time-averaged enstrophy equal to the enstrophy in the original
Navier-Stokes flow, i.e.,
\begin{equation} 
\label{eq:EE}
\left[\tE(\cdot;{\varphi})\right]_T = \Big[\E(\cdot)\Big]_T =: \E_0.
\end{equation}
This condition {can be interpreted as fixing the
  $L^2([0,T];L^2(\Omega))$ norm of the vorticity field in the LES
  flow. It} thus defines the following nonlinear manifold in the space
$\cS$, cf.~\eqref{eq:S}, of nondimensional functions parameterizing
the eddy viscosity
\begin{equation}
{\M} := \left\{ {\varphi} \in \cS \: : {\left[\tE(\cdot;{\varphi})\right]_T}=\E_0, \right\}.
\label{eq:M2}
\end{equation}

The optimal eddy viscosity $\nuc = \nuc(s)$ can then be obtained from
the solution of the following optimization problems via ansatz
\eqref{eq:nu} with the function $\varphi$ replaced with $\cvarphi$,
{where as in the previous section we consider the unconstrained
  and constrained formulation with and without the nonlinear
  constraint} {\begin{problem}\label{prob:nu} For the system
    \cref{eq:LES} and objective functional \cref{eq:J2}, find
\begin{align} 
\widecheck{{\varphi}}:=& \underset{{\varphi} \in {\cS}} {\argmin} \, \J({\varphi}), \label{prob:minJS} \tag{\showprobC{A}} \\
\widecheck{{\varphi}}:=& \underset{{\varphi} \in {\M}} {\argmin} \, \J({\varphi}). \label{prob:minJM} \tag{\showprobC{B}}
\end{align}
\end{problem}
}

We emphasize that Problem \ref{prob:nu} has a fundamentally
different structure than Problems \ref{prob:Heat1} and
\Ref{prob:Heat2}, since the control variable $\varphi(s/\smax)$ is a
function of the {\em dependent} (state) variable $s$,
cf.~\eqref{eq:nu}, rather than the {\em independent} variables ($t$ or
$x$) as is typical in PDE-constrained optimization problems. In other
words, solving Problem \ref{prob:nu} can be interpreted as finding
an optimal form of the nonlinearity in the closure model.  As
described in the next section, this aspect will have significant
ramifications for how Problem \ref{prob:nu} is solved.

\section{Solution Approach\label{sec:solution}}

In this section we first briefly present an adjoint-based approach to
compute the gradient {$\grad_{\phi} {\J}(\phi)$} of the (reduced)
objective functional in Problem \ref{prob:Heat1} and \ref{prob:Heat2},
which is quite standard, and then describe how an analogous approach
can be used to define the projection operator {$P_{\TM_{\phi}}$}
realizing the {orthogonal} projection onto the subspace tangent
to the constraint manifold ${\M}$, cf.~\eqref{eq:desc1} and Figure
\ref{fig:M}. Finally, we provide some details how these approaches can
be adapted to {perform analogous tasks in the solution of}
Problem \ref{prob:nu}. Throughout these derivations particular
attention will be paid to ensuring the required regularity of the
obtained solutions which will be done using the framework established
in \cite{pbh04}.

\subsection{Evaluation of the Cost Functional Gradients}
\label{sec:grad}

In order to determine the gradient $\grad_{\phi}{\J}$ of the objective
functional \eqref{eq:J}, we begin by computing its G\^{a}teaux
(directional) differential in the direction of some arbitrary
perturbation $\phi' \in H^1(0,T)$ of the heat flux
\begin{equation} 
\label{eq:dJ}
{\J}'(\phi; \phi') = \frac{d}{d\epsilon} \, {\J} \left(\phi + \epsilon \phi'\right) \Big|_{\epsilon = 0} = \int_{0}^T [u(\phi)|_{b} - \ub] \, u'(x, t; \phi, \phi') \, dt,
\end{equation}
where $u'(t, x; \phi, \phi')$ satisfies the system obtained as a
perturbation of the governing system \eqref{eq:Heat}
\begin{subequations} \label{eq:pert}
\begin{alignat}{3} 
{\K} u' := &\p{u'}{t} - \Delta u'& &= 0& \qquad (t, x) &\in (0, T] \times \Omega, \\
&\p{u'}{x}\Big|_{x=a}& &= \phi'(t)& \quad t &\in (0,T],  \label{eq:pertBC} \\
&\p{u'}{x}\Big|_{x=b}& &= 0& \quad t &\in (0, T], \\
&u'(t=0)& &= 0& \quad x &\in \Omega. 
\end{alignat} 
\end{subequations}

In order to extract the gradient $\grad_{\phi}{\J}$ from the G\^{a}teaux
differential \eqref{eq:dJ}, we use the fact that the differential is a
bounded linear functional on both $L^2(0,T)$ and $H^1(0,T)$ when
viewed as a function of $\phi'$, and invoke the Riesz representation
theorem to obtain \cite{b77}
\begin{equation} 
\label{eq:Riesz}
{\J}'({\phi}; \phi')  = \left\langle {\grad_{\phi}^{L^2}}{\J}, \phi' \right\rangle_{L^2{(0, T)}} 
= \Big\langle {\grad_{\phi}} {\J}, \phi' \Big\rangle_{H^1{(0, T)}},
\end{equation}
where the gradients ${\grad_{\phi}^{L^2}} {\J} \in L^2(0,T)$ and ${\grad_{\phi}} {\J} \in
H^1(0,T)$ are the Riesz representers in the two spaces.  While in
\cref{eq:desc} we require the gradient in the space $H^1(0, T)$, it
is convenient to first obtain the gradient with respect to the $L^2$
topology.  Since the expression on the RHS of \eqref{eq:dJ} is not
consistent with the Riesz form \eqref{eq:Riesz} as the perturbation
$\phi'$ does not appear explicitly in it, but is {instead} hidden in the
boundary condition \eqref{eq:pertBC}, we introduce the {\em adjoint
  field} $u^* \; : \; [0,T] \times \Omega \rightarrow \RR$ and define
the following duality-pairing relation
{
\begin{equation}
\begin{aligned}
  \left( {\K} u' , u^* \right) :&=  \int_0^T \intO {(\K u') \, u^*} \, dx \, dt  \\
&  = \int_0^T \intO {u' \,  (\K^* u^*}) \, dx \, dt -  \overbrace{\int_{0}^T [u(\phi)|_{b} - \ub]
  \, u'(b, t; \phi, \phi') \, dt}^{{{\J}}'({{\phi}}; {\phi}')} -
\int_{0}^{T} u^*\big|_{x=a} \, \phi'(t) \, dt = 0,
\end{aligned}
\label{eq:dual}
\end{equation}}
where integration by parts was performed with respect to both space
and time and the field $u^*$ solves {\em adjoint system}
\begin{subequations} \label{eq:adj}
\begin{alignat}{3} 
{\K^* u^*} := -&\p{u^*}{t} - \Delta u^*& &= 0& \qquad (t, x) &\in (0, T] \times \Omega, \\
&\p{u^*}{x}\Big|_{x=a}& &= 0& \quad t &\in (0, T], \\
&\p{u^*}{x}\Big|_{x=b}& &= u(\phi)|_{b} - \ub& \quad t &\in (0, T], \\
&u^*({t=T})& &= 0& \quad x &\in \Omega. 
\end{alignat} 
\end{subequations}
\\
Collecting \eqref{eq:pert}, \eqref{eq:dual} and \eqref{eq:adj}, we
obtain an expression for the G\^{a}teaux differential {consistent
  with the Riesz form \eqref{eq:Riesz}, namely,}
${{\J}}'({{\phi}}; {\phi}') =  - \int_{0}^{T} u^*\big|_{x=a} \, \phi'(t) \, dt$.
The gradient defined with respect to the $L^2$ topology is then deduced using the first equality in \eqref{eq:Riesz}
\begin{equation} 
\label{eq:gradL2}
\grad_{\phi}^{L^2}{{\J}}(t) =-u^*(t, x)\big|_{x=a}.
\end{equation}

The required Sobolev gradient ${\grad_{\phi} {\J}}$ is then obtained using the
second equality in \eqref{eq:Riesz} and the definition \eqref{eq:ipHp}
which gives
\begin{equation}
{\J}'(\phi; \phi') = \int_{0}^{T} {\grad_{\phi}^{L^2} {\J}} \, \phi' \, dt 
= \int_{0}^{T} {\grad_{\phi} {\J}} \, \phi' \, dt + \ell^2 \, \int_{0}^{T} \frac{d({\grad_{\phi} {\J}})}{dt}\, \frac{d\phi'}{dt} \, dt.
\label{eq:dJr}
\end{equation}
{Performing integration by parts with respect to time $t$ in the
  second term and assuming that both} the perturbation and the
gradient vanish at the endpoints of the time interval,
$\phi'(0) = \phi'(T) = 0$ and
${\grad_{\phi} {\J}}\big|_{t=0} = {\grad_{\phi} {\J}}\big|_{{t=T}} =
0$, we obtain the Sobolev gradient as the solution of the elliptic
boundary-value problem \cite{pbh04}
\begin{subequations}
\label{eq:gradH1} 
\begin{alignat}{3} 
&\left[\Id - \ell^2 \, \frac{d^2}{dt^2} \right] {\grad_{\phi} {\J}}(t) & &= {\grad_{\phi}^{L^2} {\J}}(t) \qquad t \in (0,T), \label{eq:gradH1a}  \\
&{\grad_{\phi} {\J}}\big|_{t=0} = {\grad_{\phi} {\J}}\big|_{t={T}} & &= 0.
\end{alignat} 
\end{subequations}
Given the form of the discrete gradient flow \eqref{eq:desc}, this
ensures that the behavior of the control variable at the endpoints
$t = 0$ and $t = T$ can be prescribed in the initial guess $\phi_0$
(in other words, with the Sobolev gradients defined as in
\eqref{eq:Riesz} and \eqref{eq:gradH1}, the {discrete} gradient flow
\eqref{eq:desc} does not affect the boundary values of the control
variable, which is desirable in some applications). It can be shown
\cite{pbh04} that extraction of Sobolev gradients via the
boundary-value problem \eqref{eq:gradH1} can be interpreted as
application of a low-pass filter to the $L^2$ gradient with $\ell$
acting as the cut-off parameter (it is thus a smoothing operation
where Fourier components of the gradient with wavelengths shorter than
$\ell$ are {damped}).

\subsection{Projection onto the Subspace Tangent to the Constraint Manifold \label{sec:proj}}

Projection onto the subspace ${\TM_{\phi}}$ tangent to the manifold
${\M}$ at a given element $\phi \in {\M}$ is a key step in enforcing
the associated constraint, cf.~Figure \ref{fig:M}. Since {in the
  present problem} the codimension of the manifold is {one},
cf.~\eqref{eq:M}, the tangent subspace can be characterized in terms
of an element ${\N_{\phi}} \in H^1(0,T)$ normal to it in a suitable
sense (if the codimension of the manifold is higher,
$\codim({\M}) > 1$, then there would be $\codim({\M})$ such elements
and the considerations below generalize in a natural way). We note
that the projection operator {$P_{\TM_{\phi}}$} in the discrete
gradient flow \eqref{eq:desc} needs to be defined with respect to the
$H^1$ inner product, cf.~\eqref{eq:ipHp}, as this is the topology of
the ambient space. However, as we did above, we will first obtain the
normal element {$\N_{\phi}^{L^2}$} defined with respect to the $L^2$
inner product and then deduce its Sobolev counterpart, both as
functions of time.

We begin by considering the G\^{a}teaux differential of the constraint
$\left[E(\cdot; \phi)\right]_T = E_0$, cf.~\eqref{eq:M}, which yields
\begin{equation}
\label{eq:dE0}
\left[E'(\cdot; \phi, \phi')\right]_T := \frac{d}{d\epsilon} \ \left[E(\cdot; \phi + \epsilon {\phi}')\right]_T \bigg|_{\epsilon = 0} 
= \frac{1}{T}\int_0^T \intO u(t, x; \phi) \, u'(t, x; \phi, \phi') \, dx \, dt = 0, 
\end{equation}
where $ u'(t, x; \phi, \phi')$ is the solution of the perturbation
system \eqref{eq:pert}. We then seek to express the directional
differential \eqref{eq:dE0}, which is a bounded linear functional of
$\phi'$ on $L^2(0,T)$ and $H^1(0,T)$, in terms of Riesz identities in
these spaces, i.e.,
\begin{equation} 
\label{eq:Riesz2}
\left[E'(\cdot; \phi, \phi')\right]_T
= \left\langle {\N_{\phi}^{L^2}}, {\phi}' \right\rangle_{L^2{(0, T)}} = \Big\langle {\N_{\phi}}, {\phi}' \Big\rangle_{H^1{(0, T)}} = 0.
\end{equation}
Noting that the expression on the RHS of \eqref{eq:dE0} does not
explicitly depend on the perturbation $\phi'$ {which instead
  appears in the boundary condition \eqref{eq:pertBC}}, we introduce a
{\em new} adjoint state
$v^* \; : \; [0,T] \times \Omega \rightarrow \RR$ that will satisfy
the following adjoint system
\begin{subequations} \label{eq:adj2}
\begin{alignat}{3} 
{\K}^* v^* := -&\p{v^*}{t} - \Delta v^*& &= u(t, x; \phi)& \qquad (t, x) &\in (0, T] \times \Omega, \\
&\p{v^*}{x}\Big|_{x=a}& &= 0& \quad t &\in (0, T], \\
&\p{v^*}{x}\Big|_{x=b}& &= 0& \quad t &\in (0, T], \\
&v^*({t=T})& &= 0& \quad x &\in \Omega. 
\end{alignat} 
\end{subequations}
Then, {introducing} the duality relation 
\begin{equation}
  {
\left( {\K}u', v^* \right)
:=  \int_0^T \intO ({\K}u')  v^* \, dx dt  = \overbrace{\int_0^T \intO u' ({\K^* v^*}) \, dx \, dt}^{\left[E'(\cdot; \phi, \phi')\right]_T} + \int_{0}^{T} v^*\big|_{x=a} \, \phi'(t) \, dt = 0,}
\label{eq:dualE0}
\end{equation}
we conclude that
$\left[E'(\cdot; \phi, \phi')\right]_T = - \int_{0}^{T} v^*\big|_{x=a}
\, \phi'(t) \, dt$. Using the first equality in \eqref{eq:Riesz2}
{provides an expression for the element normal to the tangent
  subspace $\TM_{{\phi}}$ with respect to the $L^2$ inner product,
  namely,}
\begin{equation} 
\label{eq:gradL2E0}
{\N_{\phi}^{L^2}}(t) = -v^*(t, a), \qquad t \in [0,T].
\end{equation}
Using the second equality in \eqref{eq:Riesz2} {and following the
  same steps as above}, we obtain ${\N_\phi} \in H^1(0,T)$, {the
  element normal to $\TM_{{\phi}}$ with respect to the $H^1$
  topology,} as the solution of the elliptic boundary-value problem
\eqref{eq:gradH1} with the source term on the RHS in
\eqref{eq:gradH1a} replaced with expression \eqref{eq:gradL2E0}.
{For a given $\phi \in \M$,} the projection operator then takes
the form
\begin{equation}
\label{eq:PTM}
{\forall \ {\Phi} \in H^1(0,T) \qquad 
{P_{\TM_{\phi}} \Phi} := {\Phi - \zeta \, \N_{\phi}}, \quad \text{where} \ \
{\zeta := \frac{\left\langle \Phi, \N_{\phi} \right\rangle_{H^1(0, T)}}{\left\langle \N_{\phi}, \N_{\phi} \right\rangle_{H^1(0, T)}} }}
\end{equation}
can be interpreted as the Lagrange multiplier {associated with
  the condition 
  $\left\langle P_{\TM_{\phi}} \Phi, \N_{\phi} \right\rangle_{H^1(0,T)} = 0$.} We
emphasize that each iteration of the gradient descent \eqref{eq:desc}
requires the solution of two adjoint systems: the solution of system
\eqref{eq:adj} allows us to determine the gradient of the objective
functional, whereas the solution of \eqref{eq:adj2} is needed in order
to construct the projection operator \eqref{eq:P}. The two adjoint
systems are defined in terms of the same operator {$\K^*$}, but are
subject to different boundary conditions and source terms.

As already indicated in {\Cref{sec:Introduction}}, for general
nonlinear manifolds ${\M}$ the projected discrete gradient flow
\eqref{eq:desc} produces local minimizers which do not satisfy the
constraint exactly, but only with an error of order
$\mathcal{O}(\tau_n^2)$ at each iteration. If the constraint is
homogeneous as is the case in Problem {\ref{prob:Heat1}}, then
one can define the retraction operator
$\R_{\M} \; : \: {\TM} \rightarrow {\M}$ in terms of simple
normalization
\begin{equation}
\label{eq:retract}
\R_{{\M}}(\phi) := \sqrt{\frac{E_0}{\left[E(\cdot; \phi)\right]_T}} \ \phi.
\end{equation}
The resulting Riemannian gradient flow then satisfies the constraint
exactly (i.e., up to round-off errors only) {\cite{ams08}}. {The
  steps required to solve \ref{prob:Heat1M} and
  \ref{prob:Heat2M} with the approach described above are summarized
  as Algorithm \ref{algo1}. As is evident from this algorithm, the
  additional per-iteration computational cost resulting from
  approximate enforcement of the constraint
  $\left[E(\cdot; \phi)\right]_T = E_0$ consists in the solution of
  the adjoint system \eqref{eq:adj2} and the boundary-value problem
  \eqref{eq:gradH1}, lines 9--12.}

\begin{algorithm}[H] 
{
  \caption{
{
Solution of \ref{prob:Heat1M} and
    \ref{prob:Heat2M}, cf.~Figure \ref{fig:M}
  \newline \textbf{Input:} \\
  \hspace*{0.3cm} \makebox[1.1cm]{$\ub(t)$\hfill} --- prescribed profile \\
  \hspace*{0.3cm} \makebox[1.1cm]{$E_0$\hfill} --- value of the constraint \\
  \hspace*{0.3cm} \makebox[1.1cm]{$\Delta x,\Delta t$\hfill} --- numerical discretization parameters \\
  \hspace*{0.3cm} \makebox[1.1cm]{$\ell$\hfill} --- Sobolev length scale\\
  \hspace*{0.3cm} \makebox[1.1cm]{$\epsilon_{\J}$\hfill} --- tolerance in the termination criterion \\
  \hspace*{0.3cm} \makebox[1.1cm]{$\phi_0(t)$\hfill} --- initial guess \\
  \textbf{Output:} \\
  \hspace*{0.3cm} \makebox[1.1cm]{$\cphi(t)$\hfill} --- optimal flux \\ 
}
}
\label{algo1}
\begin{algorithmic}[1]
\State \textbullet~set $n = 0$
\State \textbullet~set $\phi^{(0)} = \phi_0$
\Repeat 
\State \textbullet~set $n = n+1$
\State \textbullet~solve the governing system \eqref{eq:Heat} with the
boundary condition \eqref{eq:HeatBC} given by $\phi^{(n)}$
\State \textbullet~solve the adjoint system \eqref{eq:adj}
\State \textbullet~determine the {$L^2$} gradient $\grad_{\phi}^{L^2}\J(\phi^{(n)})$, cf.~\eqref{eq:gradL2}
\State \textbullet~determine the Sobolev gradient $\grad_{\phi}\J(\phi^{(n)})$, cf.~\eqref{eq:gradH1}
\State \textbullet~solve the adjoint system \eqref{eq:adj2}
\State \textbullet~determine the normal element in {$L^2$} $\N_{\phi}^{L^2}$, cf.~\eqref{eq:gradL2E0}
\State \textbullet~determine the normal element in the Sobolev space $H^1$ $\N_{\phi}$, cf.~\eqref{eq:gradH1}
\State \textbullet~determine the projection of the Sobolev gradient onto the tangent subspace $\TM_{\phi^{(n)}}$, \eqref{eq:PTM}
\State \textbullet~determine the optimal step length $\tau^{(n)}$ by
solving a line-minimization problem \\ \quad\qquad while applying retraction
\eqref{eq:retract} (\ref{prob:Heat1M} only)
\State \textbullet~update the flux $\phi^{(n)} \ \leftarrow
\R_{\M}\left( \phi^{(n)} - \tau_n \, {P_{\TM_{\phi^{(n)}}} \phi}^{(n)}  \right)$
\Until $\left(\J(\phi^{(n-1)}) - \J(\phi^{(n)}) \right) /
\J(\phi^{(n-1)}) < \epsilon_{\J}$
\end{algorithmic}
}
\end{algorithm}

\subsection{Solution of \Cref{prob:nu}}
\label{sec:nu}

Here we provide some details describing how the framework introduced
in Sections \ref{sec:grad} and \ref{sec:proj} can be used to solve
\Cref{prob:nu}. {Since the system considered in this problem is
  more complicated,} to balance completeness and brevity, we only
state the most important steps and refer the reader to
\cite{Matharu2022a} for {further} details.

{First, we compute the gradient of the objective functional
  \eqref{eq:J2} and} following the same procedure as stated in
\Cref{sec:grad} begin by {determining} the G\^{a}teaux
differential of \cref{eq:J2}
\begin{equation} 
\label{eq:dJ2}
{\J'({\varphi}; {\varphi}') := \frac{d}{d\epsilon} \ \J\Big({{\varphi}} + \epsilon {\varphi}'\Big) \Big|_{\epsilon = 0}
= \frac{1}{D}\,\int_0^T \intO \left( \P(t) - \tP(t; {\varphi}) \right) \bfDelta \tw(t, \x; {\varphi}) \, \tw'(t, \x; {\varphi}, {\varphi}') \, d\x \, dt,}
\end{equation}
where ${\varphi}'$ is an arbitrary perturbation of
${\varphi} \in {\cS}$, and $\tw'(t, \x; {\varphi}, {\varphi}')$
satisfies the corresponding perturbation system {obtained by
  linearizing the governing system \eqref{eq:LES}, namely,}
\begin{subequations}
	\label{eq:Pert}
	\begin{align}
	&\K\begin{bmatrix} \tw' \\  \\ \tpsi' \end{bmatrix} := 
	 \begin{bmatrix}
	\partial_t \tw' + \gradperp \tpsi' \cdot \gtw + \gradperp \tpsi \cdot \gtw' + \a \tw' \qquad \qquad \qquad \\ 
	- \bfgrad \cdot \left(2 (\gtw \cdot \gtw') (\frac{d \nu}{d{s}} \, \varphi \gtw + \frac{\eta^3 \, \sqrt{s} + \nu_0}{\smax}\, \frac{d \varphi}{d\sigma} \, {\gtw} )+ (\nuN + \nu) \gtw' \right)\\
	\bfDelta \tpsi' + \tw' 
	\end{bmatrix} \\
	&\phantom{\K\begin{bmatrix} \tw' \\  \\ \tpsi' \end{bmatrix} :}= \begin{bmatrix} \bfgrad \cdot \left( \, (\eta^3 \, \sqrt{s} + \nu_0) {\varphi}'  \gtw \, \right) \\  \\ 0\end{bmatrix}, \label{eq:Pert_eqn} \\
	&\tw'({t=0, \x}) = 0, \label{eq:Pert_IC}
	\end{align}
\end{subequations}
where for ease of notation, we have denoted $\sigma:= s/\smax$. Defining the adjoint system as 
\begin{subequations}
	\label{eq:AdjJ2}
	\begin{align}
	&\K^*\begin{bmatrix} \tw^* \\  \\ \tpsi^* \end{bmatrix} := 
	 \begin{bmatrix}
	-\partial_t \tw^* - \gradperp \tpsi \cdot \gtw^* + \a \tw^* + \tpsi^* \qquad \qquad \qquad \qquad \\ 
	- \bfgrad \cdot \left(2 \, (\gtw \cdot \gtw^*) \, (\frac{d \nu}{d{s}} \, \varphi \gtw + \frac{\eta^3 \, \sqrt{s} + \nu_0}{\smax}\, \frac{d \varphi}{d\sigma} \, {\gtw} ) + (\nuN + \nu) \gtw^*\right)\\ 
	\bfDelta \tpsi^* - \gradperp \cdot (\tw^* \, \gtw)
	\end{bmatrix} = \begin{bmatrix}  W \\  \\ 0\end{bmatrix}, \label{eq:Adjeqn} \\
	&\tw^*({t=T, \x}) = 0, \label{eq:AdjIC}
	\end{align}
\end{subequations}
with source term $W(t,\x) := \frac{1}{D}\,\left( \P(t) - \tP(t; {\varphi}) \right) \bfDelta \tw(t, \x; {\varphi})$, we obtain the duality-pairing relation
\begin{equation}
\begin{aligned}
\left( \K\begin{bmatrix} \tw' \\ \tpsi' \end{bmatrix}, \begin{bmatrix} \tw^* \\ \psi^* \end{bmatrix} \right)
:= & \int_0^T \intO \K\begin{bmatrix} \tw' \\ \tpsi' \end{bmatrix} \cdot \begin{bmatrix} \tw^* \\ \tpsi^* \end{bmatrix} \, d\x \, dt \\
=& \int_0^T \intO \begin{bmatrix} \tw' \\ \tpsi' \end{bmatrix} \cdot \K^* \begin{bmatrix} \tw^* \\ \tpsi^* \end{bmatrix} \, d\x \, dt = \left( \begin{bmatrix} \tw' \\ \tpsi' \end{bmatrix}, \K^*\begin{bmatrix} \tw^* \\ \psi^* \end{bmatrix} \right),
\end{aligned}
\label{eq:dual2}
\end{equation}
where integration by parts was performed with respect to both space
and time {and all boundary terms vanish due to periodic boundary
  conditions. This relation together with the adjoint system
  \eqref{eq:AdjJ2} allow us to re-express the G\^{a}teaux differential
  \cref{eq:dJ2} as}
\begin{align} 
\J'({\varphi}; {\varphi}') &= \phantom{-}\int_0^T \intO W(t,\x) \,
                             \tw'(t, \x; {\varphi}, {\varphi}') \, d\x \, dt,    \nonumber\\
&= -\int_0^T \intO \, \left(\eta^3 \, \sqrt{s} + \nu_0 \right) \, \left(\gtw \cdot \gtw^*\right)
\, {\varphi}' \, d\x \, dt, \nonumber \\
&= \int_0^1 \left[-\int_0^T \intO \, \delta\left(\frac{\gtw \cdot \gtw}{\smax} - \sigma \right) \, \left(\eta^3 \, \sqrt{s} + \nu_0 \right) \,\gtw \cdot \gtw^* \, d\x \, dt \right] \,  \varphi'(\sigma) \, d\sigma,
\label{eq:dJ2b}
\end{align}
where the substitution
{$\varphi'(\gtw \cdot \gtw) = \int_0^1 \, \delta\left(\frac{\gtw \cdot
      \gtw}{\smax} - \sigma \right) \, \varphi'(\sigma) \, d\sigma$}
{with $\delta(\cdot)$ denoting the Dirac delta distribution was
  made} and Fubini's theorem was used to swap the order of
integration. {These last two steps are needed in order to change
  the integrations variables in \eqref{eq:dJ2b} from $d\x \, dt$ to
  $ds$, such that the G\^{a}teaux differential \cref{eq:dJ2} can be
  expressed in the required Riesz form}
\begin{align} 
\label{eq:RieszJ2}
\J'(\varphi; \varphi') = \Big\langle \grad_{\varphi}\J, \varphi' \Big\rangle_{H^2({[0, 1]})} &= \left\langle \grad_{\varphi}^{L^2}\J, {\varphi}' \right\rangle_{L^2({[0, 1]})}.
\end{align}
{This then allows} us to extract the gradient with respect to the
$L^2$ topology as
\begin{equation} 
\label{eq:gradL2b}
\grad_{\varphi}^{L^2}{\J}(\sigma) 
=-\int_0^T \intO \, \delta\left(\frac{\gtw \cdot \gtw}{\smax} - \sigma
\right) \, \left(\eta^3 \, \sqrt{s} + \nu_0 \right) \, \gtw \cdot
\gtw^* \, d\x \, dt, \qquad {\sigma \in [0,1].}
\end{equation}

As was done in \Cref{sec:grad}, we use the {Riesz
  representation in the Sobolev space $H^2(0,T)$,
  cf.~\eqref{eq:Riesz2}, and the definition of the subspace $\cS$ in
  \eqref{eq:S} which after integration by parts with respect to
  $\sigma$, cf.~\eqref{eq:dJr},} give the Sobolev gradient
$\grad_{\varphi}\J$ as the solution of the following elliptic
boundary-value problem
\begin{subequations}
	\label{eq:gradH2BVP}
	\begin{align}
	\left[\Id - \ell_1^2 \, \frac{d^2}{d\sigma^2} + \ell_2^4 \frac{d^4}{d\sigma^4}\right] \grad_{\varphi}\J(\sigma) &= \grad_{\varphi}^{L^2}\J(\sigma), \qquad \sigma \in [0, 1], \label{eq:gradH2} \\
	\frac{d^{(1)} \, (\grad_{\varphi}\J)}{d\sigma^{(1)}} \Big|_{\sigma=0,1} &= \frac{d^{(3)} \, (\grad_{\varphi}\J)}{d\sigma^{(3)}} \Big|_{\sigma=0,1} = 0. \label{eq:gradH2bc}
	\end{align}
\end{subequations}
{We add that the boundary conditions in \eqref{eq:gradH2bc}, see
  also \eqref{eq:S}, were selected in order to ensure the vanishing
  of the terms at $\sigma = 0,1$ which result from integration by
  parts.}

Similarly, we determine the normal element $\N_{\varphi}$ by
considering the G\^{a}teaux differential of the constraint
\cref{eq:EE} with respect to $\varphi$ and invoking the Riesz
representation theorem, i.e.,
\begin{align} 
\left[\tE'(\cdot;{\varphi}, {\varphi}')\right]_T := \frac{d}{d\epsilon} \ \left[\tE(\cdot;{\varphi} + \epsilon {\varphi}')\right]_T \Bigg|_{\epsilon = 0} &= \frac{1}{T} \int_0^T \intO \tw(t, \x; {\varphi}) \, \tw'(t, \x; {\varphi}, {\varphi}') \, d\x \, dt, \nonumber \\
&= \Big\langle \N_{\varphi}, \varphi' \Big\rangle_{H^2({[0, 1]})} = \left\langle \N_{\varphi}^{L^2}\J, {\varphi}' \right\rangle_{L^2({[0, 1]})} = 0. \label{eq:RieszE0}
\end{align}
Introducing the new {\em adjoint fields} $\tbw^*$ and $\tbpsi^*$,
which satisfy the adjoint system \cref{eq:AdjJ2}, with the source term
$W(t,\x) := \frac{1}{T}\,\tw(t, \x; {\varphi})$, we {can} use the
duality pairing \cref{eq:dual2} to conclude that
\begin{align*} 
\left[\tE'(\cdot;{\varphi}, {\varphi}')\right]_T &= -\int_0^T \intO \, \left(\eta^3 \, \sqrt{s} + \nu_0 \right) \, \left(\gtw \cdot \gtbw^*\right)
\, {\varphi}' \, d\x \, dt,\\
&= \int_0^1 \left[-\int_0^T \intO \, \delta\left(\frac{\gtw \cdot \gtw}{\smax} - \sigma \right) \, \left(\eta^3 \, \sqrt{s} + \nu_0 \right) \,\gtw \cdot \gtbw^* \, d\x \, dt \right] \,  \varphi'(\sigma) \, d\sigma.
\end{align*}
{Then,} using the second equality in \cref{eq:RieszE0}, we obtain
{an expression for the element normal to the subspace
  $\TM_{\varphi}$ with respect to the $L^2$ topology as}
\begin{equation} 
\label{eq:gradL2E0b}
\N_{\varphi}^{L^2}(\sigma) =-\int_0^T \intO \, \delta\left(\frac{\gtw
    \cdot \gtw}{\smax} - \sigma \right) \, \left(\eta^3 \, \sqrt{s} +
  \nu_0 \right) \, \gtw \cdot \gtbw^* \, d\x \, dt, \qquad  {\sigma \in [0,1].}
\end{equation}
{Finally, employing the last equality in \eqref{eq:RieszE0} and
  performing analogous steps as described above, we obtain the normal
  element $\N_\varphi \in H^2(0,1)$ as a solution of the
  boundary-value problem \cref{eq:gradH2BVP}, but with the source term
  in \eqref{eq:gradH2} replaced with the RHS of
  \cref{eq:gradL2E0b}. This then allows us to construct the projection
  operator $P_{\TM_{\varphi}} $ in the form \eqref{eq:PTM}.} {The
  steps required to solve \Cref{prob:nu} are analogous to Algorithm
  \ref{algo1}, with obvious modifications.}

\section{Numerical Discretization and Validation}
\label{sec:numer}

In this section we briefly describe the numerical discretizations
employed to solve Problems \ref{prob:Heat1}, \ref{prob:Heat2} and
\ref{prob:nu} using the discrete gradient flow \eqref{eq:desc}. Then
we validate key elements of these computations requiring solution of
adjoint systems, namely, evaluation of the cost functional gradients
\eqref{eq:gradL2} and \eqref{eq:gradL2b} and the normal elements
\eqref{eq:gradL2E0} and \eqref{eq:gradL2E0b}. {The proposed approach has been implemented in MATLAB and the code is available in \citep{MatharuGitHeatOptConstr} (it can be used to generate the figures shown in \Cref{sec:numerP12,sec:resultsP12}).}

\subsection{Problems \ref{prob:Heat1} and \ref{prob:Heat2}}
\label{sec:numerP12}

The governing system \eqref{eq:Heat} and the adjoint systems
\eqref{eq:adj} and \eqref{eq:adj2} are discretized in space with a
standard second-order finite-difference scheme and integrated in time
using a second-order Crank-Nicolson scheme.

\begin{figure}
\centering
\begin{tikzpicture}[line cap=round,line join=round,>=triangle 45]
\begin{axis}[
x=0.8cm,y=0.8cm,
xmin=0,
xmax=15,
ymin=0,
ymax=10,
ticks=none,
axis lines=left,
x label style={at={(axis description cs:.5,-.01)},anchor=north},
y label style={at={(axis description cs:-.01,.5)},anchor=south},
xlabel={\large $\log_{10}\epsilon$},
ylabel={\large $\log_{10}|1-\kappa(\epsilon)|$}]
\clip(-4.4237396135265685,-5.212224881355937) rectangle (16.539571761208386,10);
\draw [<-, line width=2pt,color=black!40!green] (1,7) -- (4,4);
\draw [-, loosely dashed, line width=2pt,color=black!40!green] (1.5,6.5)-- (7,1) node [midway, below, rotate=-45] {\large Round-Off Errors, $\mathcal{O}(\epsilon^{-1})$};
\draw [line width=2pt,color=red] (4,4)-- (10,4) node [midway, below] {\large Gradient Errors};
\draw [-, loosely dashed, line width=2pt,color=blue] (7,1) -- (12.5,6.5) node [midway, below, rotate=45] {\large Truncation Errors $\mathcal{O}(\epsilon^1)$};
\draw [->, line width=2pt,color=blue] (10,4)-- (13,7);
\draw [->, line width=2pt,color=blue] (10,4)-- (12.3,8);
\draw [->, line width=2pt,color=blue] (10,4)-- (11.5,9) node [midway, above, rotate=73.3] {\large $\mathcal{O}(\epsilon^p)$, $p>1$};
\draw [->,line width=2pt,color=red] (7,8.5) -- (7,4.5) node [midway, above, rotate=-90] {\large $\dt,\dx \to 0$};
\end{axis}
\end{tikzpicture}
\caption{{Schematic of the interplay between the three types of
    errors in the validation of accuracy of the gradients and of the
    normal elements based on formulas \cref{eq:kappa12,eq:kappa34},
    cf.~Section \ref{sec:numerP12}.  The slope of the right branch
    (blue) depends on the order of the finite-difference formula used
    to approximate the G\^{a}teaux differential in the denominator
    ($p = 1$ in \cref{eq:kappa12,eq:kappa34}).  }}
\label{fig:kscheme}
\end{figure}
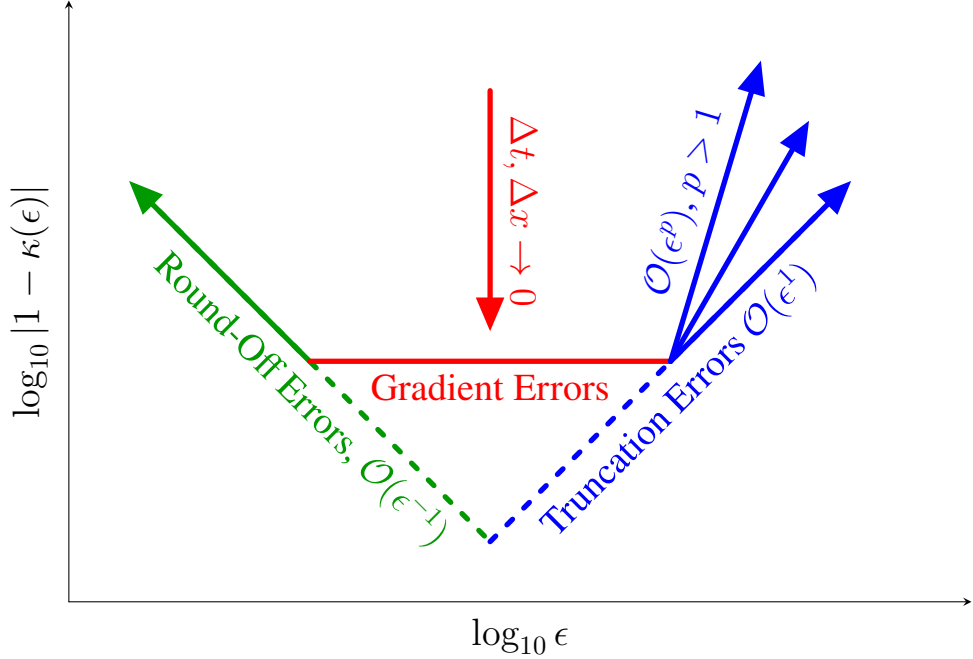

In order to validate the computation of the gradient \eqref{eq:gradL2}
and of the normal element \eqref{eq:gradL2E0} we define the following
quantities
\begin{subequations}
\label{eq:kappa12}
\begin{align}
\kappa_1(\epsilon) & := \frac{\epsilon^{-1} \left[ {\J}({\phi} + \epsilon \phi') - {\J}({\phi}) \right]}{\Big\langle {\grad_{\phi}^{L^2}} {\J}(\phi), \phi' \Big\rangle_{L^2(0, T)}}, \qquad \epsilon > 0
\label{eq:kappa1} \\
\kappa_2(\epsilon) & := \frac{\epsilon^{-1} \left[ E(\cdot; {\phi} + \epsilon \phi') - E(\cdot; {\phi}) \right]_T}{\Big\langle {\N_{\phi}^{L^2}}, \phi' \Big\rangle_{L^2(0, T)}}{,}
\label{eq:kappa2}
\end{align}
\end{subequations}
in which the numerators represent first-order finite-difference
approximations of the G\^{a}teaux differentials \eqref{eq:dJ} and
\eqref{eq:dE0}, whereas the denominators are the corresponding Riesz
forms \eqref{eq:Riesz} and \eqref{eq:Riesz2}. Clearly, we expect that
$\kappa_i(\epsilon) \approx 1$, $i=1,2$, and any deviations of these
quantities from unity {results from a combination of errors of
  the following three distinct types:
  \begin{itemize}
  \item errors in the numerical solution of the PDE systems
    \eqref{eq:Heat}, \eqref{eq:adj} and \eqref{eq:adj2}, and in the
    evaluation of the integrals in \eqref{eq:M} and \eqref{eq:J};
    these errors result from the discretization of \revtt{these} equations
    and expressions in space and in time; and are therefore controlled
    by the corresponding discretization parameters ($\Delta x$ and
    $\Delta t$, but in general there may also be some other numerical
    parameters) and are independent of $\epsilon$; hereafter, we refer
    to these errors as ``gradient errors'',
\item
truncation errors in the finite-difference formula {in
  \eqref{eq:kappa1}--\eqref{eq:kappa2}} which are
proportional to $\epsilon$ (if a finite-difference formula of order $p
> 1$ were used in \eqref{eq:kappa1}--\eqref{eq:kappa2}, then the
truncation errors would be proportional to $\epsilon^p$), and 
\item
subtractive cancellation (round-off) errors proportional to $\epsilon^{-1}$.
\end{itemize}
An interplay of these types of errors in formulas such as
\eqref{eq:kappa1}--\eqref{eq:kappa2} is shown schematically in
\Cref{fig:kscheme}. Evidently, both the truncation and subtractive
cancellation (round-off) errors are artefacts of the structure of
these formulas and accuracy of the gradient $\grad_{\phi}^{L^2} {\J}$
and of the normal element $\N_{\phi}^{L^2}$ is controlled solely by
the gradient errors. Therefore, in order to validate their
computation, one must show that these errors vanish as the numerical
parameters are refined, i.e., as $\Delta x, \Delta t \rightarrow 0$,
which is manifested by the lowering of the plateau labelled ``Gradient
Errors'' in \Cref{fig:kscheme} until it ultimately disappears. We add
that the order of the finite-difference formula used to approximate
the G\^{a}teaux differentials in the denominators of
\eqref{eq:kappa1}--\eqref{eq:kappa2} does not affect the accuracy with
which the gradient $\grad_{\phi}^{L^2} {\J}$ and the normal element
$\N_{\phi}^{L^2}$ are computed (it only affects the slope of the
branch corresponding to large $\epsilon$ marked in blue in
\Cref{fig:kscheme}).}

{To fix attention, we focus here on the more general Problem
  \ref{prob:Heat2} and add that very similar results were also
  obtained for Problem \ref{prob:Heat1}. In our tests we set} $a = 0$,
$b = 1$, $T = 1$ and
\begin{equation}
u_0(x) = 10 \, \cos(\pi \, x), \quad x \in [0,1], \qquad\qquad   
{\bar{\phi}(t) = e^{-4}\,t\,\left[ 18 - 1000\,\cos\left(\frac{15\pi}{2}\, t\right) \right]}, \quad t \in [0,T],
\label{eq:Heattrue}
\end{equation}
where $\bar{\phi}(t)$ is the ``true'' flux defining, via the solution
of system \eqref{eq:Heat}, the target profile $\bar{u}_b(t)$ appearing
in the error functional \eqref{eq:J}. In the diagnostic quantities
\eqref{eq:kappa1}--\eqref{eq:kappa2} we use
\begin{equation}
{\phi_0 := \phi(t) = 18 \, \sin\left( \frac{\pi}{2}\, t \right) e^{-4t}}, \qquad
\phi'(t) = 4t \, e^{(-4 + \pi)t}, \qquad t \in [0,T],
\label{eq:Heatphi}
\end{equation}
representing, respectively, the ``point'' (in the function space
$H^1(0,T)$) where the G\^{a}teaux differentials are computed and the
direction.

\begin{figure}
  \centering
  \subfigure[]
  {
    \includegraphics[scale=0.6]{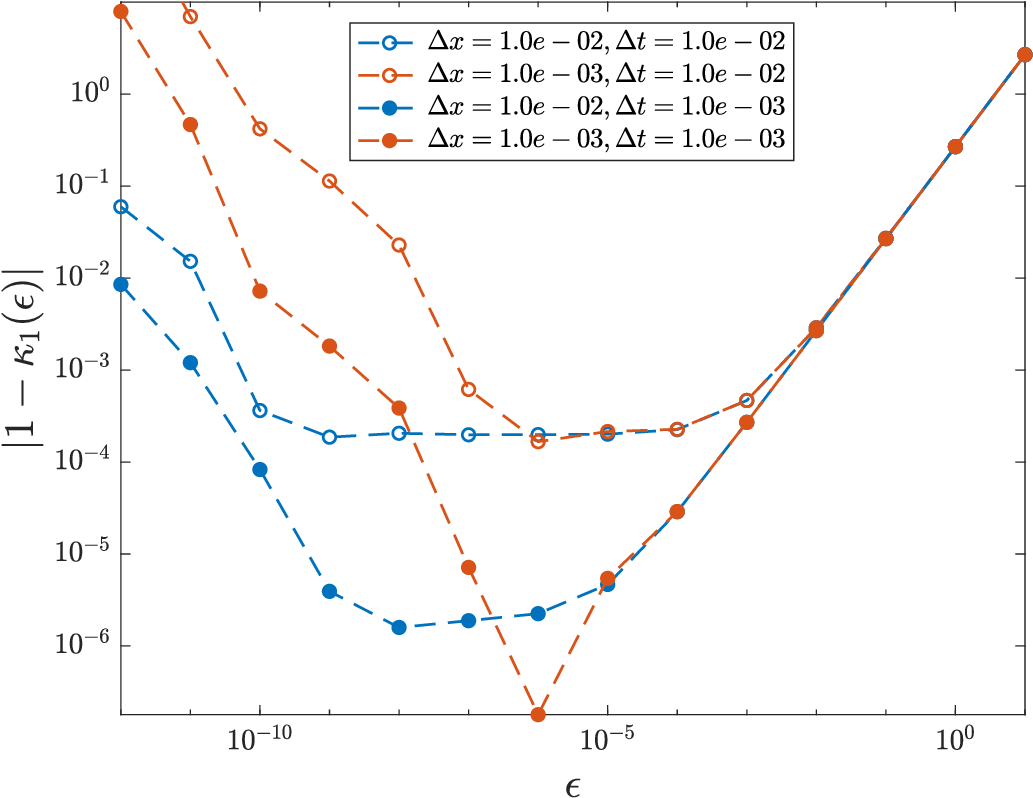}
    \label{fig:Kap1log}
  }
  \subfigure[]
  {
    \includegraphics[scale=0.6]{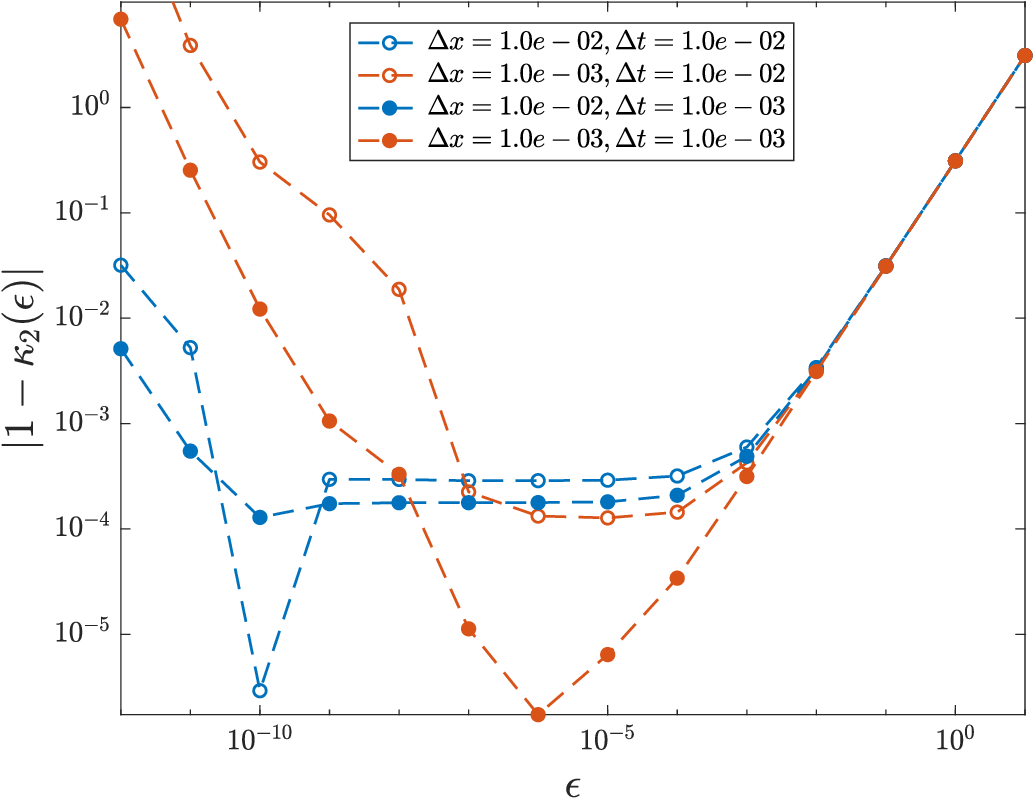}
    \label{fig:Kap2log}
  }
  \caption{[Problem \ref{prob:Heat2}]
    {Dependence of \subref{fig:Kap1log} $|1 -
    {\kappa_1(\epsilon)}|$ and \subref{fig:Kap2log} $|1 -
    {\kappa_2(\epsilon)}|$ on $\epsilon$ for different spatial and
    temporal discretizations $\Delta x$ and $\Delta t$.} }
\label{fig:kappa12}
\end{figure}
In Figures \ref{fig:kappa12}a and \ref{fig:kappa12}b we show the
dependence of the quantities {$|1 - {\kappa_i(\epsilon)}|$}, $i=1,2$,
on $\epsilon$ for different indicated values of the space and time
discretization parameters $\Delta x$ and $\Delta t$, {which
  reveals the expected behavior, cf.~\Cref{fig:kscheme}. In
  particular, we} observe that $\kappa_1(\epsilon)$ and
$\kappa_2(\epsilon)$ deviate from the unity for very small and very
large values of $\epsilon$ which is due to, respectively, the
subtractive cancellation (round-off) errors and the truncation errors
in the finite-difference formula {in
  \eqref{eq:kappa1}--\eqref{eq:kappa2}}. However, we also observe that
for intermediate values of $\epsilon$ {spanning several orders of
  magnitude} both $\kappa_1(\epsilon)$ and $\kappa_2(\epsilon)$
exhibit plateaus {corresponding to gradient errors of order
  $\mathcal{O}\left(C_{\Delta x}(\Delta x)^2 + C_{\Delta t} (\Delta
    t)^2 \right)$ for some $C_{\Delta x},C_{\Delta t} > 0$ (the
  structure of this error reflects the order of accuracy of the
  discretization techniques used to approximate relation
  \eqref{eq:Heat}, \eqref{eq:M}, \eqref{eq:J}, \eqref{eq:adj} and
  \eqref{eq:adj2}). Inspection of the results in Figure
  \ref{fig:kappa12}a shows that the gradient errors are reduced
  resulting in the lowering of the plateaus when $\Delta t$ is
  refined, indicating that the accuracy of gradient evaluations is
  controlled by the time discretization (i.e., the gradient error is
  dominated by the term proportional to $(\Delta t)^2$ and errors
  resulting from the discretization in space are relatively
  unimportant). On the other hand, the results in Figure
  \ref{fig:kappa12}b show that the accuracy with which the normal
  element is evaluated is more or less equally controlled by both
  $\Delta x$ and $\Delta t$.  In both cases, the plateaus ultimately
  disappear when the resolution is sufficiently refined indicating
  that the gradient errors become smaller than the subtractive
  cancellation (round-off) and truncation errors.  This} demonstrates
the convergence of the approximations of the cost functional gradient
and the normal element. {Since the Sobolev gradients and normal
  directions are obtained by solving the elliptic boundary value
  problem \cref{eq:gradH1} and can be {therefore viewed as
    low-pass filtered versions of the corresponding} elements found in
  $L^2$ \cite{pbh04}, {their correctness follows from the
    validation provided above together with relations \eqref{eq:Riesz}
    and \eqref{eq:Riesz2}.  In}} the computations presented in Section
\ref{sec:resultsP12} the Sobolev gradients and normal directions are
computed using the inner product \eqref{eq:ipHp} with
{$\ell_1 = 0.01$}.


\subsection{\Cref{prob:nu}}
\label{sec:numerP3}

The Navier-Stokes system \eqref{eq:2DNS}, the {corresponding} LES
system \eqref{eq:LES} and the {adjoint systems \eqref{eq:AdjJ2}} are
approximated using a standard Fourier pseudo-spectral method in space
and using a third-order IMEX scheme introduced in \cite{Alimo2020} in
time. Evaluation of the closure terms involving the state-dependent
eddy viscosity \eqref{eq:nu} requires interpolation from the physical
space $\Omega$ to the state space $\I$ and differentiation with
respect to {the state variable} $s$, steps which are performed
using methods based on Chebyshev polynomials
\cite{trefethen2013approximation}. We refer the reader to
\cite{Matharu2022a} for further numerical details.

The ``target'' Navier-Stokes flow $w$ is defined as the solution of
system \eqref{eq:2DNS} obtained with the parameters
{$\nu_N = 4 \times 10^{-4}$, $\alpha = 5 \times 10^{-3}$, $F = 2$,
  $k_{a} = k_{b} = 4$ over the time window $T = 50 \approx 27.8\, t_e$},
where
$t_e := \left[ \int_{0}^{T} \E(t) \, dt / ( 8\pi^2 T)\right]^{-1/2}$
is the eddy turnover time \citep{Bracco2010}. {This time window
  is chosen to make \Cref{prob:nu} physically interesting, but at the
  same time is not too long as to make the computation of gradients
  problematic due to exponential divergence of nearby trajectories (it
  is worth mentioning here that reliable computation of sensitivities
  of long-time averages of quantities defined for chaotic systems is
  facilitated by the ``shadowing approach''
  \cite{WangHuBlonigan2014,BloniganWang2018}).}  {The direct numerical
  simulation (DNS) of \cref{eq:2DNS} and the LES computation solving
  \cref{eq:LES} use, respectively, $N_x = 256$ and $N_x = 64$
  equispaced grid points in each spatial direction. They also use
  $N_s = 128$ Chebyshev points to discretize} the state space $\I$.
The LES system \eqref{eq:LES} is obtained using a ``box'' filter with
the cutoff {$k_c = 8$} \cite{Matharu2022a}. The use of this rather
``aggressive'' filter is dictated by the desire to define a problem
where the presence of the constraint \eqref{eq:EE} has a significant
effect on the solutions of Problem \ref{prob:nu}, which will in turn
allow us to more easily elucidate the properties of the proposed
approach.

In order to validate the computation of the cost functional gradient
and the normal element, we define the following diagnostic quantities
analogous to \eqref{eq:kappa1}--\eqref{eq:kappa2}
\begin{subequations}
\label{eq:kappa34}
\begin{align}
\kappa_3(\epsilon) & := \frac{\epsilon^{-1} \left[ \J({\varphi} + \epsilon \varphi') - \J({\varphi}) \right]}{\Big\langle {\grad_{\varphi}^{L^2}\J(\varphi)}, \varphi' \Big\rangle_{L^2(0, 1)}}, \qquad \epsilon > 0 \label{eq:kappa3} \\
\kappa_4(\epsilon) & := \frac{\epsilon^{-1} \left[\tE(\cdot;{\varphi} + \epsilon {\varphi}') - \tE(\cdot;{\varphi})\right]_T}
{\Big\langle {\N^{L^2}_{{\varphi}}}, {\varphi}' \Big\rangle_{L^2(0, 1)}},
\label{eq:kappa4}
\end{align}
\end{subequations}
where $\varphi = \varphi(s)$ and $\varphi' = \varphi'(s)$ are the
nondimensional functions corresponding, via ansatz \eqref{eq:nu}, to
the eddy viscosity in the Leith model {\citep{Leith1968,
    Leith1971, Leith1996}}
\begin{equation}
\label{eq:nu0}
\nu_{L}(s) = (C_{L} \, k_{c})^3 \, \sqrt{s}
\end{equation}   
with the constant $C_{L}=4.702 \times 10^{-3}$ chosen to ensure that
$\varphi \in {\M}$, and to the perturbations
\begin{subequations}
\begin{align}
\nu'(s) &= \nu_{L}(s), \qquad \qquad \qquad \text{with} \quad C_{L}=4.2 \times 10^{-3}, \label{Constr:pert1} \\
\nu'(s) &= \exp\left(\frac{-2s}{30}\right). \label{Constr:pert3}
\end{align}
\end{subequations}

\begin{figure}
  \centering
  \subfigure[]
  {
    \includegraphics[scale=0.6]{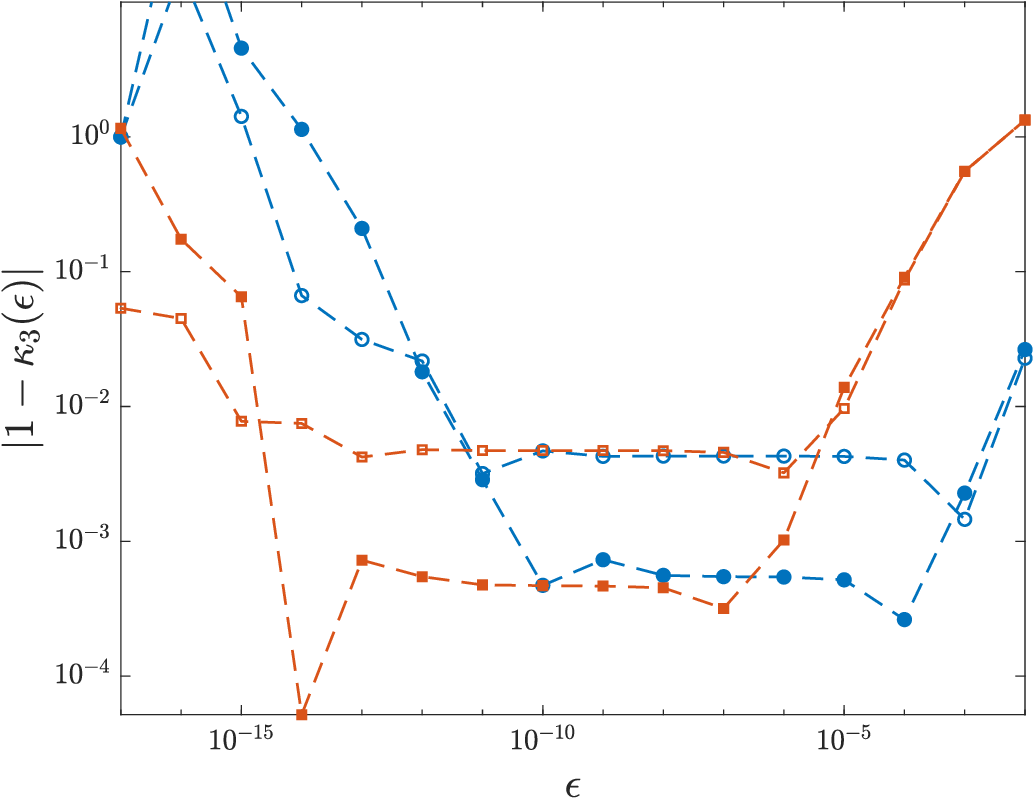}
    \label{fig:Kap2Palin_log}
  }
  \subfigure[]
  {
    \includegraphics[scale=0.6]{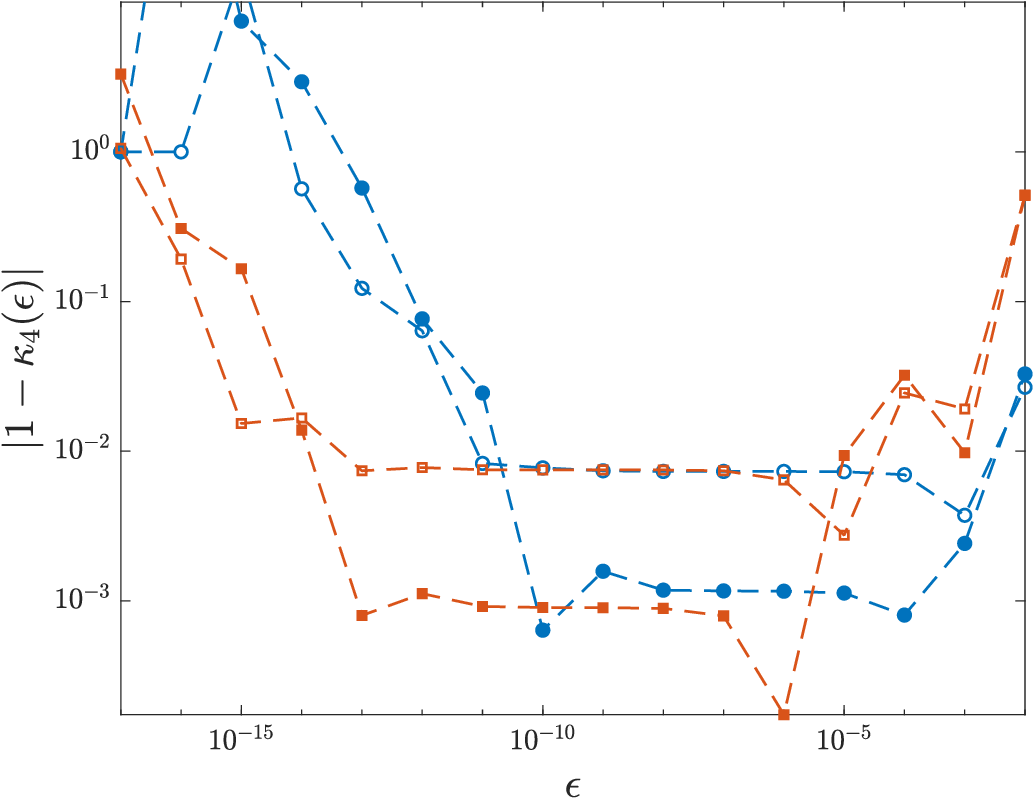}
    \label{fig:Kap3Const_log}
  }
  \caption{[Problem \ref{prob:nu}:] {Dependence of (a)
      $|1 - {\kappa_3(\epsilon)}|$ and (b)
      $|1 - {\kappa_4(\epsilon)}|$ on $\epsilon$ for different
      temporal discretizations with $\dt = 10^{-2}$ (empty
      symbols) and {$\dt = 10^{-3}$} (filled symbols). The
      results are shown for the perturbations $\varphi'$ corresponding
      to $\nu'$ given in {\cref{Constr:pert1} (blue circles)
        and \cref{Constr:pert3} (red squares).} }}
\label{fig:kappa34}
\end{figure}

The quantities {$|1-\kappa_3(\epsilon)|$ and $|1-\kappa_4(\epsilon)|$}
are shown as functions of $\epsilon$ for different perturbations
$\varphi'$ and two different time steps $\dt$ in
Figures~\ref{fig:kappa34}a and \ref{fig:kappa34}b, respectively. Since
the spatial discretization {(with respect to $\x$) and the
  discretization of the state space $s$ we use are both spectrally
  accurate {(i.e., they have} exponentially small errors)}, we focus
here on assessing the effect of the time discretization, which has an
algebraic only rate of convergence, on the accuracy of the gradient
and the normal element.  In Figures \ref{fig:kappa34}a and
\ref{fig:kappa34}b we see the expected behavior of the quantities
$\kappa_3(\epsilon)$ and $\kappa_4(\epsilon)$,
{cf.~\Cref{fig:kscheme}. In particular, they approach unity for
  intermediate values of $\epsilon$ as the time step $\dt$ is refined
  indicating the vanishing of the gradient error which is dominated by
  the time discretization. However, in contrast to the results of
  analogous tests for Problem \ref{prob:Heat2}, cf.~Figure
  \ref{fig:kappa12}a,b, now the convergence is slower with the
  gradient error of order $\mathcal{O}(\Delta t)$, which is a
  consequence of the difficulties in evaluating the integrals
  involving the Dirac distributions in \eqref{eq:gradL2b} and
  \eqref{eq:gradL2E0b}.  We can nevertheless conclude that the
  approximations of the cost functional gradient and of the normal
  element are convergent.}  In the computations presented in Section
\ref{sec:resultsP3} the Sobolev gradients and normal directions are
computed using the inner product \eqref{eq:ipHp} with {$\ell_1 = 1000$
  and $\ell_2 = 100$}.

\section{Results\label{sec:results}}

In this section we study solutions to Problems \ref{prob:Heat1},
\ref{prob:Heat2} and \ref{prob:nu} focusing our attention on how well
and how efficiently the constraints involved are satisfied. In
iterations \eqref{eq:desc} the step size $\tau_n$ is determined using
Brent's algorithm to solve the line-minimization problem
\cite{pftv86}.  In the unconstrained versions of the problems we use
the Polak-Ribi{\`{e}}re version of the conjugate gradient method
\cite{nw00} to accelerate convergence of iterations. In the
constrained versions of the problems {when no retraction operator is available}, we limit the magnitude of the
step size $\tau_n$ to reduce drift {away from the constraint
  manifold $\M$}.

\subsection{Problems \ref{prob:Heat1} and \ref{prob:Heat2}}
\label{sec:resultsP12}

In Problems \ref{prob:Heat1} and \ref{prob:Heat2} we attempt to
reconstruct the flux $\phi(t)$ at the left boundary of the domain
{by minimizing the error functional \eqref{eq:J}} such that the
temperature at the right boundary matches the prescribed target
temperature $\bar{u}_b$ over the time window $[0,T]$; the latter
temperature is obtained by solving problem \eqref{eq:Heat} with a
certain ``true'' flux $\bar{\phi}(t)$, given in \eqref{eq:Heattrue},
applied at the left boundary. {The true flux $\bar{\phi}(t)$ thus
  represents the ``exact'' solution of Problems \ref{prob:Heat1}.A and
  \ref{prob:Heat2}.A in the sense that $\J(\bar{\phi}) = 0$ in these
  cases. However, since it does not satisfy the nonlinear
  constraint, $\bar{\phi} \notin \M$, the true flux cannot be a
  solution of Problems \ref{prob:Heat1}.B and \ref{prob:Heat2}.B.  In
  the solution of Problem \ref{prob:Heat1}.B the nonlinear constraint,
  cf.~\eqref{eq:M}, is enforced exactly using the retraction operator
  $\R_\M$ whereas in the solution of Problem \ref{prob:Heat2}.B this
  has to be done approximately using the orthogonal projection on
  the tangent subspace \eqref{eq:PTM}.}

\begin{figure}\centering
  \subfigure[]
  {
    \includegraphics[scale=0.6]{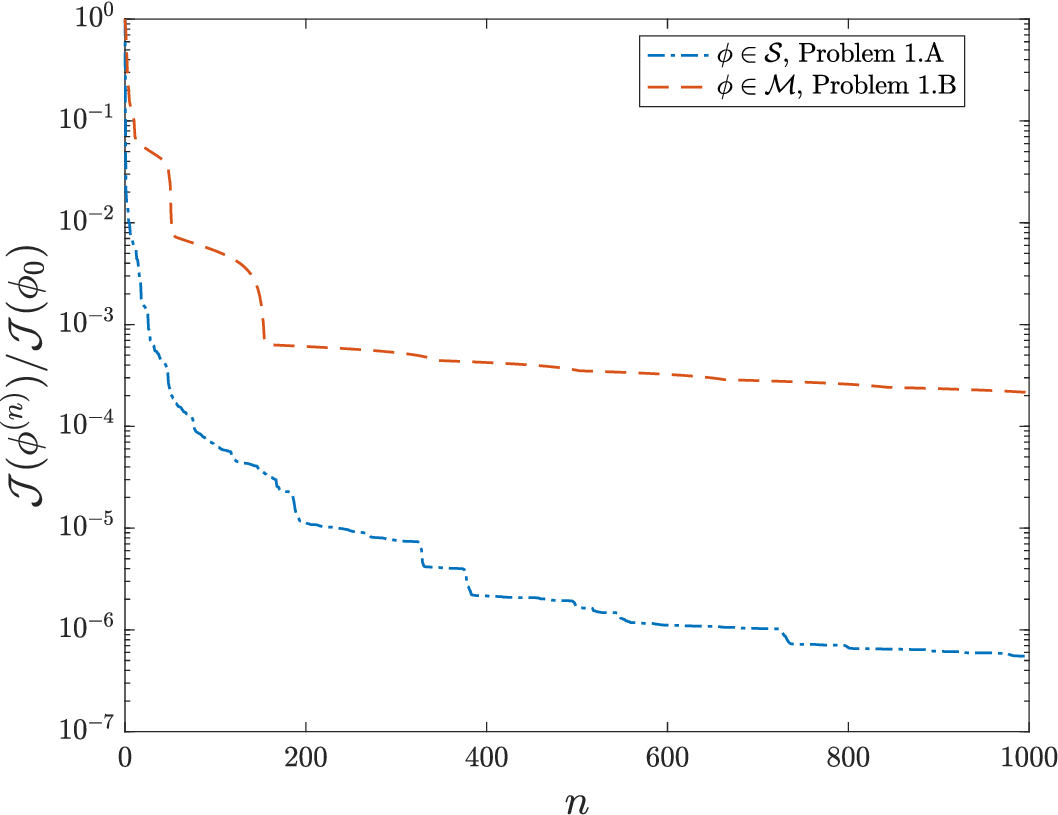}
    \label{fig:J1}
  }\quad
  \subfigure[]
  {
    \includegraphics[scale=0.6]{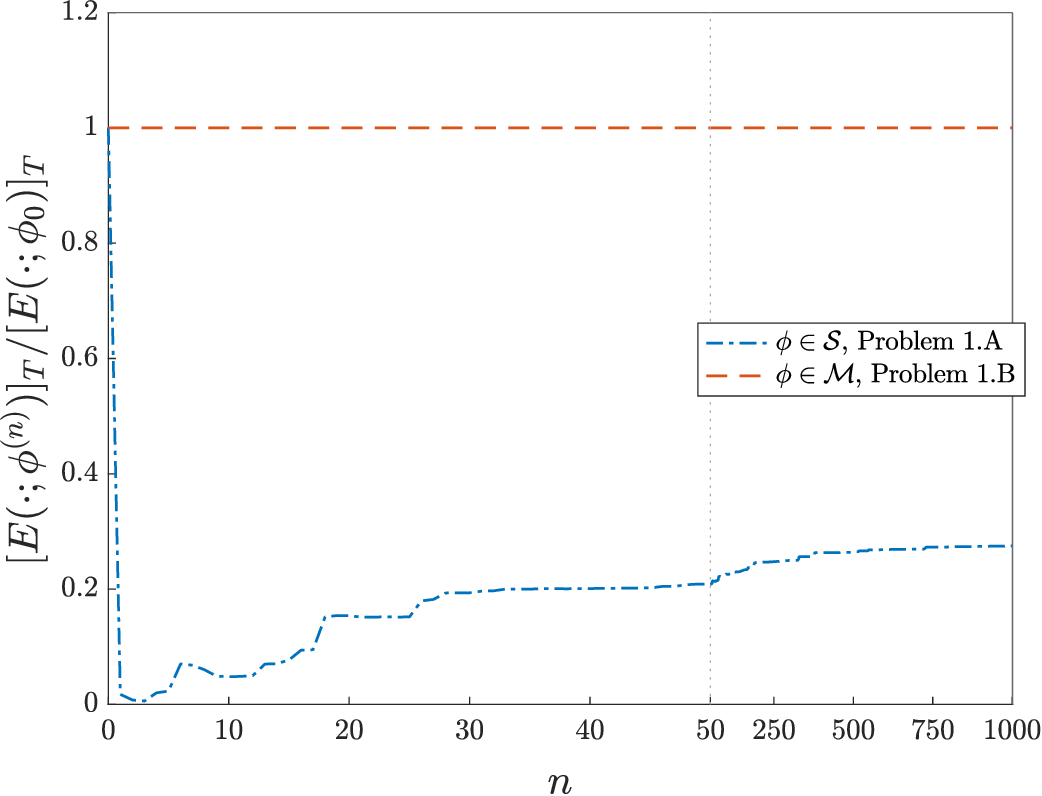}
    \label{fig:E1}
  }
  \caption{[Problem \ref{prob:Heat1}:] Dependence of (a) the
    normalized error functional
    ${{\J}}({{\phi}}^{(n)})/{{\J}}({{\phi}}_0)$ and (b) normalized
    constraint $\left[{E}\left(\cdot;
        {{\phi}}^{(n)}\right)\right]_T/\left[{E}\left(\cdot;
        {{\phi}}_0\right)\right]_T$ on the iteration $n$ in the
    solutions of the {unconstrained and constrained} optimization
    problem. {For clarity, in \subref{fig:E1}, the horizontal
      axis is split into two regions with different linear scaling
      separated by the vertical dotted line.}  }
  \label{fig:JE1}
\end{figure}

\begin{figure}\centering
    \includegraphics[scale=0.6]{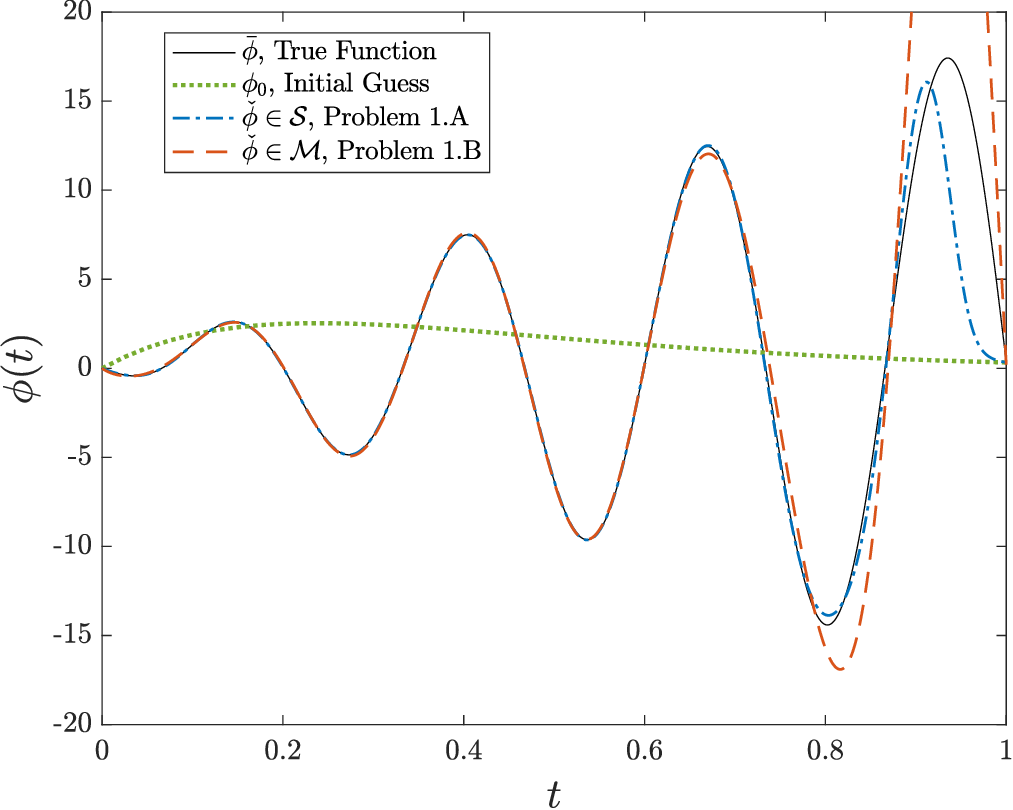}
    \caption{[Problem \ref{prob:Heat1}:] (black,
      solid line) the true flux $\bar{\phi}(t)$,
      cf.~\eqref{eq:Heattrue}, (green, dotted line) the initial guess
      $\phi_0(t)$ in \eqref{eq:desc} and \cref{eq:Heatphi}, (blue, dot-dashed line) the
      optimal flux $\cphi(t)$ in the  {unconstrained} problem, and (red, dot-solid line) the optimal flux
      $\cphi(t)$ in the {constrained} problem, all as
      functions of time $t \in [0,T]$.}
  \label{fig:cphi1}
\end{figure}

The dependence of the normalized error functional ${\J}(\phi^{(n)}) /
{\J}(\phi_0)$ and of the normalized constraint $\left[E(\cdot;
  \phi^{(n)})\right]_T / E_0$ on the iteration count $n$ in the
solution of {Problems \ref{prob:Heat1}.A and \ref{prob:Heat1}.B} is
shown, respectively, in Figures \ref{fig:J1} and \ref{fig:E1}.  The
{corresponding} optimal fluxes $\cphi(t)$ are shown as functions of
time $t \in [0,T]$ in Figure \ref{fig:cphi1}. Figure \ref{fig:J1}
shows that, as expected, in the solution of Problem \ref{prob:Heat1}.A
the error functional \eqref{eq:J} is reduced to a very low level and
this error is significantly higher in the solution of Problem
\ref{prob:Heat1}.B where an exact solution is not expected to exist.
On the other hand, the nonlinear constraint in \eqref{eq:M} is
satisfied exactly in Problem \ref{prob:Heat1}.B, cf.~Figure
\ref{fig:E1}. In contrast, in the unconstrained problem we see that
the difference $\left[E(\cdot; \phi)\right]_T - E_0$ becomes large,
{on the order of $E_0$,} after a few iterations, which is not be
surprising as this constraint is not imposed in \ref{prob:Heat1S}. As
is evident from Figure \ref{fig:cphi1}, the true flux is reconstructed
rather well in \ref{prob:Heat1S}, except for the final times where a
large oscillation is present. The oscillatory character of the optimal
flux $\cphi(t)$ in this case indicates that in its original form
Problem \ref{prob:Heat1} is ill-posed. Methods for dealing with this
issue are well developed and are usually based on Tikhonov
regularization \cite{ehn96} or formulation of the problem on spaces of
sufficiently regular functions \cite{pbh04}. Since this issue is only
tangential to the main topic of the present study, we do not consider
it further here.

\begin{figure}\centering
  \subfigure[]
  {
    \includegraphics[scale=0.6]{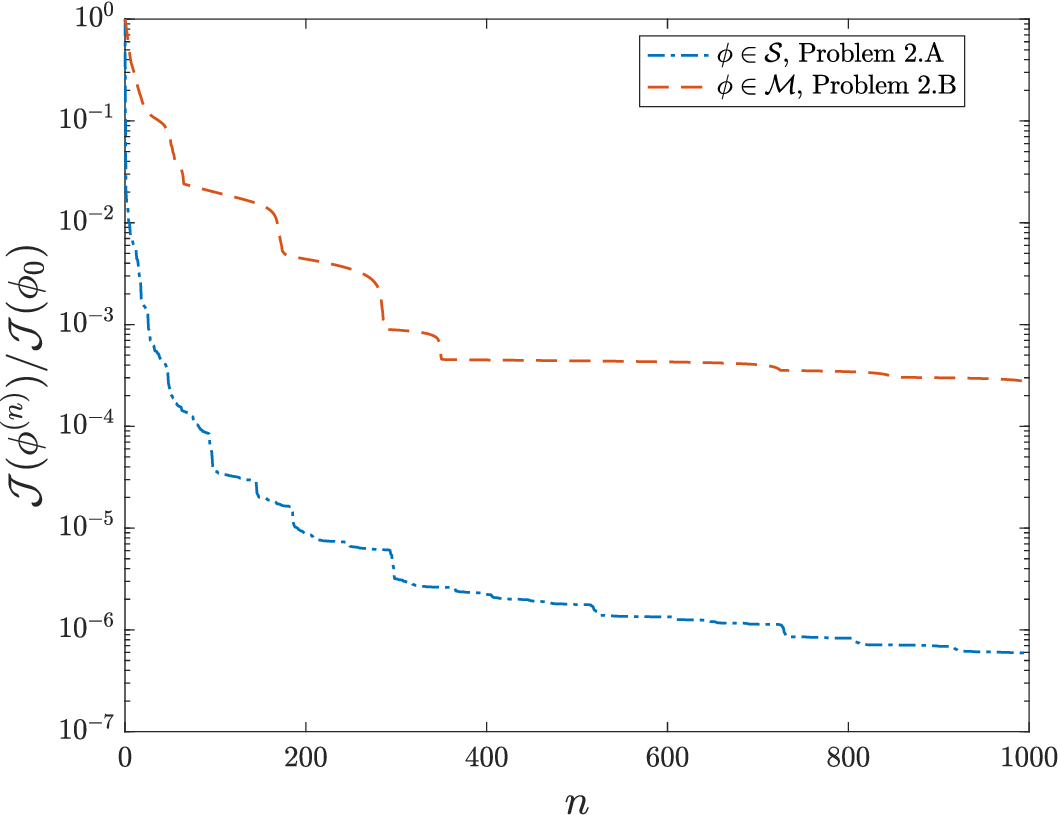}
    \label{fig:J2}
  }\quad
  \subfigure[]
  {
    \includegraphics[scale=0.6]{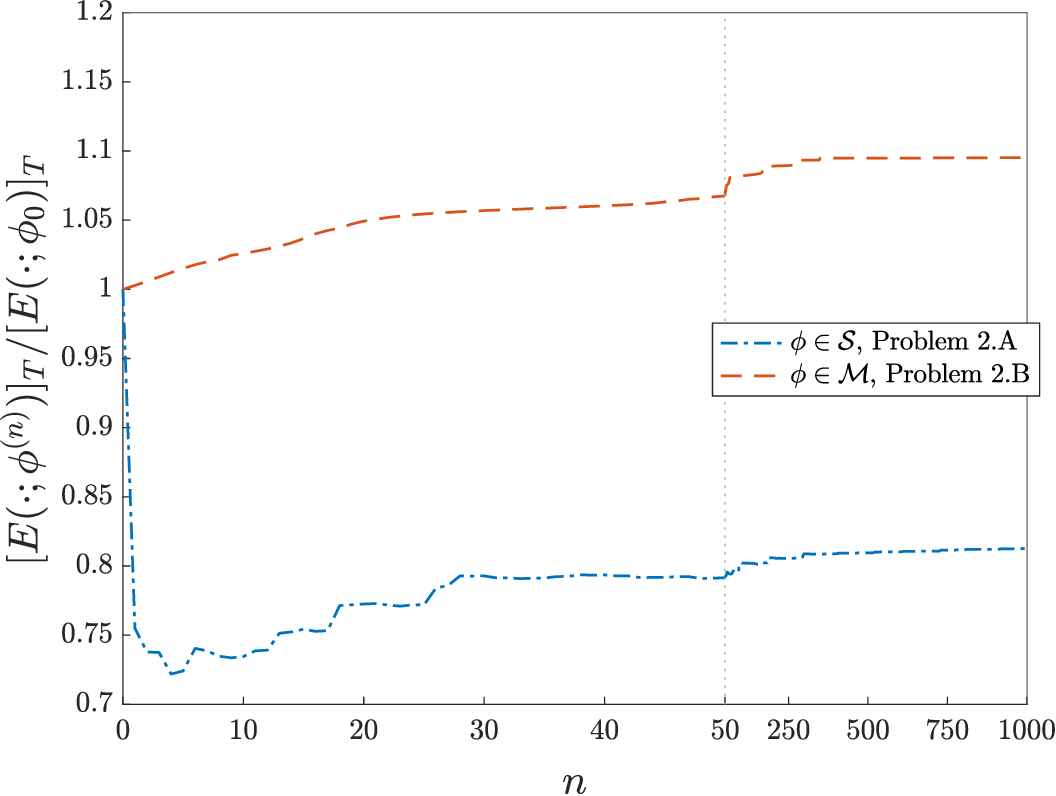}
    \label{fig:E2}
  }
  \caption{[Problem \ref{prob:Heat2}:] Dependence of (a) the
    normalized error functional
    ${{\J}}({{\phi}}^{(n)})/{{\J}}({{\phi}}_0)$ and (b) normalized
    constraint $\left[{E}\left(\cdot;
        {{\phi}}^{(n)}\right)\right]_T/\left[{E}\left(\cdot;
        {{\phi}}_0\right)\right]_T$ on the iteration $n$ in the
    solutions of the {unconstrained and constrained} optimization
    problem. {For clarity, in \subref{fig:E2}, the horizontal
      axis is split into two regions with different linear scaling
      separated by the vertical dotted line.}  }
  \label{fig:JE2}
\end{figure}

\begin{figure}\centering
    \includegraphics[scale=0.6]{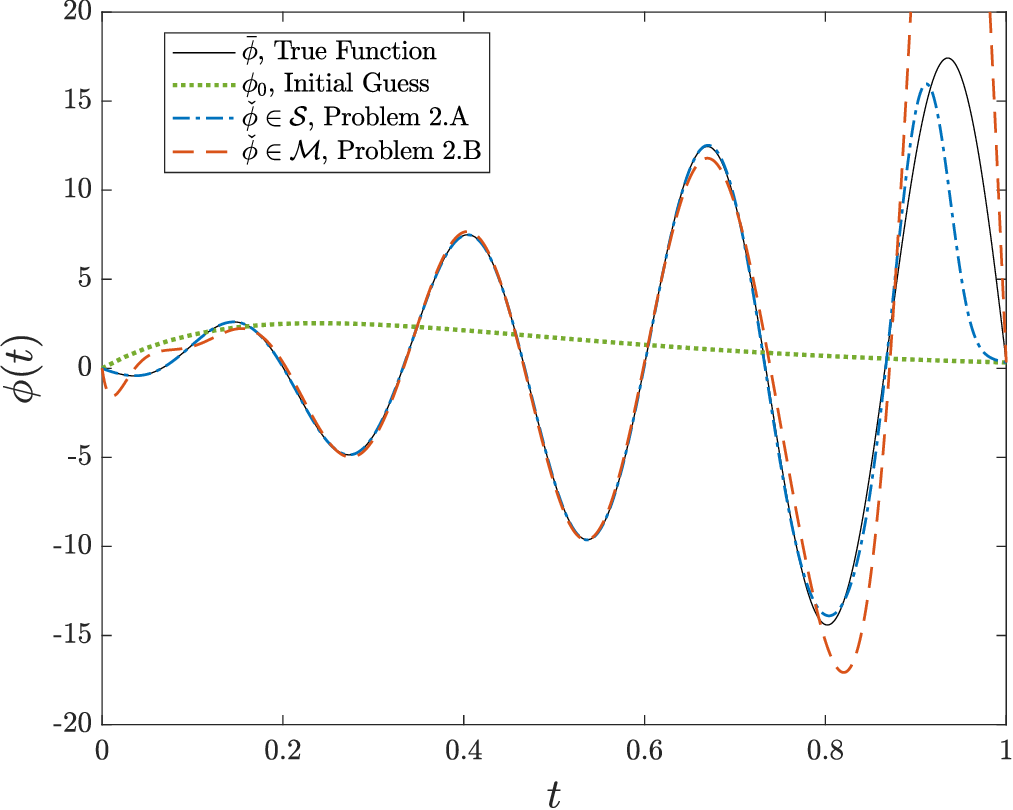}
    \caption{[Problem \ref{prob:Heat2}:] (black,
      solid line) the true flux $\bar{\phi}(t)$,
      cf.~\eqref{eq:Heattrue}, (green, dotted line) the initial guess
      $\phi_0(t)$ in \eqref{eq:desc}{ and \cref{eq:Heatphi}}, (blue, dot-dashed line) the
      optimal flux $\cphi(t)$ in the  {unconstrained} problem, and (red, dot-solid line) the optimal flux
      $\cphi(t)$ in the {constrained} problem, all as
      functions of time $t \in [0,T]$.}
  \label{fig:cphi2}
\end{figure}

{Moving on to discuss Problem \ref{prob:Heat2},} the dependence of the
normalized error functional ${\J}(\phi^{(n)}) / {\J}(\phi_0)$ and of
the normalized constraint $\left[E(\cdot; \phi^{(n)})\right]_T / E_0$
on the iteration count $n$ in shown in Figures \ref{fig:J2} and
\ref{fig:E2}.  {We emphasize that now the retraction operator $\R_\M$
  is no longer available and the nonlinear constraint in
  \ref{prob:Heat2M} can only be enforced approximately via projection
  on the tangent subspace $\TM_{\phi}$, cf.~Figure \ref{fig:M}. We can
  see that, as a result, the constraint is no longer satisfied
  exactly, but the {drift away from the constraint manifold $\M$
    is slow, much less rapid} than when the constraint is not enforced
  at all as in \ref{prob:Heat2S}. In this problem as well we observe a
  significant reduction of the error functional \eqref{eq:J},
  cf.~Figure \ref{fig:J2}, leading to a good reconstruction of the
  true flux $\bar{\phi}$ in \ref{prob:Heat2S}, cf.~Figure
  \ref{fig:cphi2}.}

\FloatBarrier

\subsection{\Cref{prob:nu}}
\label{sec:resultsP3}

\begin{figure}\centering
  \subfigure[]
  {
    \includegraphics[scale=0.6]{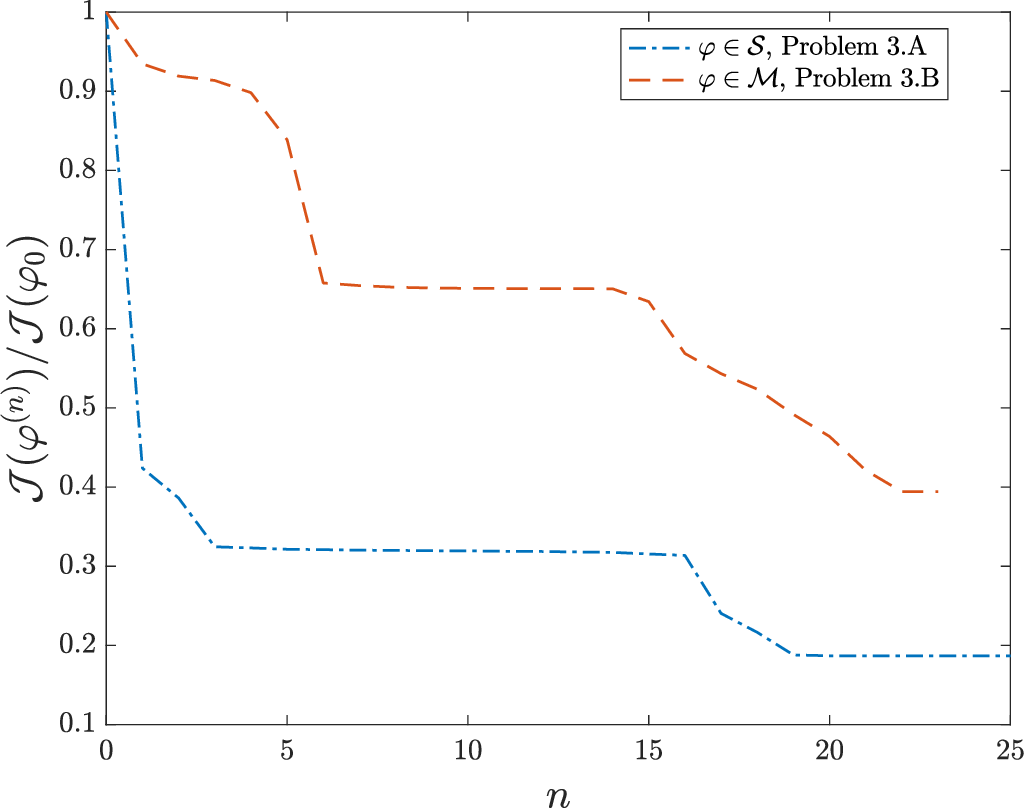}
    \label{figConstraint:costfun}
  }
  \subfigure[]
  {
    \includegraphics[scale=0.6]{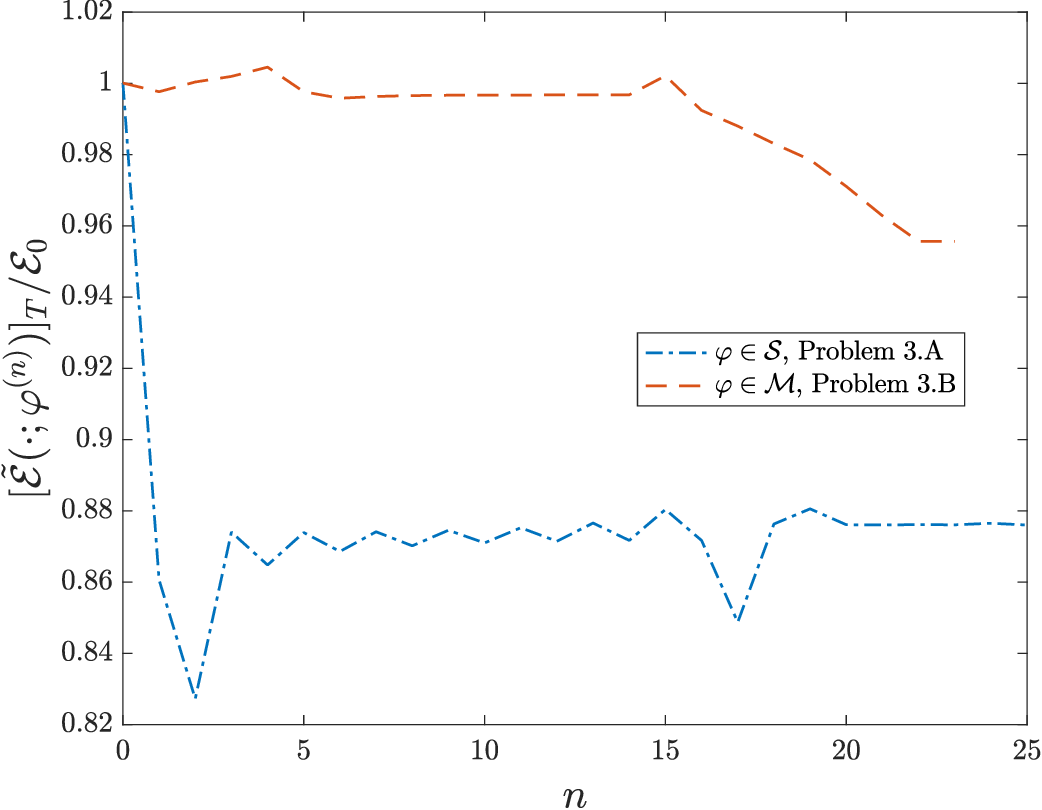}
    \label{figConstraint:EnstTW}
  }
  \caption{[Problem \ref{prob:nu}:] Dependence of (a) the normalized
    error functional ${\J}({{\varphi}}^{(n)})/{\J}({{\varphi}}_0)$ and
    (b) the normalized constraint $\left[\tE(\cdot; \varphi^{(n)})\right]_T
    / \E_0$ on the iteration $n$ for the solutions of the
    {unconstrained and constrained} optimization problem. }
  \label{fig:JE3}
\end{figure}

Problem \ref{prob:nu} is similar to Problem \ref{prob:Heat2} in that
the retraction operator is not available. {Solutions to Problems
  \ref{prob:nu}.A and \ref{prob:nu}.B} are obtained with iterations
\eqref{eq:desc} where {the eddy viscosity corresponding to the
  Leith model} \eqref{eq:nu0} is used as the initial guess.  The
dependence of the normalized error functional
$\J(\varphi^{(n)}) / \J(\varphi_0)$ and of the normalized constraint
$\left[\tE(\cdot; \varphi^{(n)})\right]_T / \E_0$ on the iteration
count $n$ is shown, respectively, in Figures \ref{fig:JE3}a and
\ref{fig:JE3}b. {In contrast to Problems \ref{prob:Heat1} and
  \ref{prob:Heat2}, Problem \ref{prob:nu} does not admit an ``exact''
  solution and as a result the values of the objective functional
  \eqref{eq:J2} attained at the minima are not very small, especially
  in the constrained Problem \ref{prob:nu}.B.  On the other hand, in
  analogy with the results for Problem \ref{prob:Heat2}, the
  minimizers found by solving the constrained Problem \ref{prob:nu}.B
  reveal a much slower drift away from the constraint manifold $\M$
  than is evident in the solution of the unconstrained Problem
  \ref{prob:nu}.A.} {Comparison of the results obtained for
  Problems \ref{prob:nu}.A and \ref{prob:nu}.B shown in Figures
  \ref{fig:JE3}a and \ref{fig:JE3}b indicates that in the particularly
  challenging (due to the constraint) latter problem a significant
  reduction of the objective functional can be achieved only if the
  iterates eventually deviate from the constraint manifold $\M$.}

\begin{figure}\centering
    \includegraphics[scale=0.6]{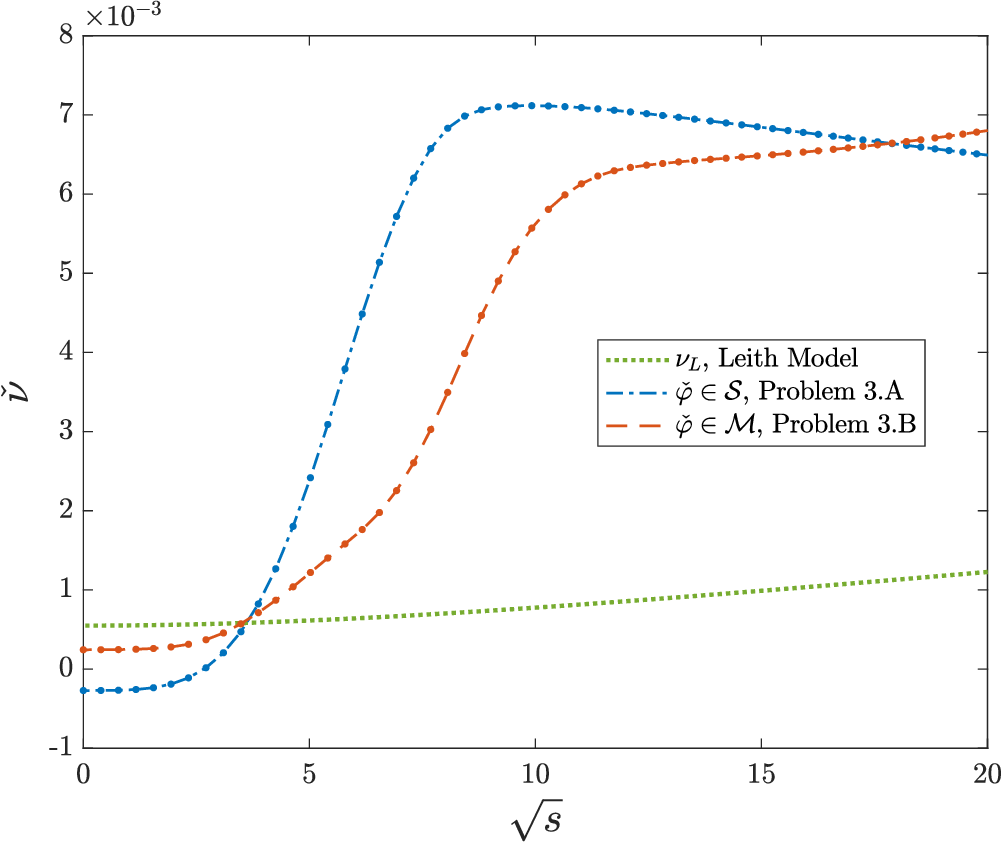}
    \caption{[Problem \ref{prob:nu}:] Dependence of the optimal eddy
      viscosity $\nuc$ on $\sqrt{s}$ for (blue, dashed-dot line)
      ${\varphi} \in \cS$ {found by solving Problem
        \ref{prob:nu}.A} and (red, dashed line) ${\varphi} \in \M$
      {found by solving Problem \ref{prob:nu}.B}.  The initial
      guess {in \eqref{eq:desc} is} given by the Leith model
      $\nu_L(s)$, cf.~\cref{eq:nu0}, {and} is shown as green
      dotted line}.
  \label{fig:nuc}
\end{figure}

\begin{figure}
  \centering
\mbox{
  \subfigure[]
  {
    \includegraphics[scale=0.5]{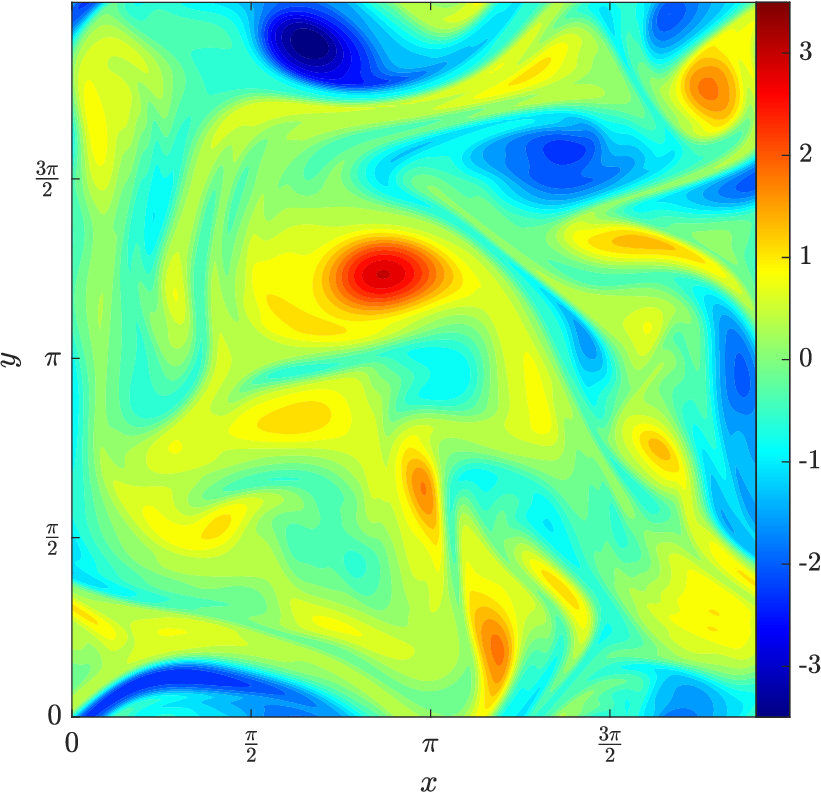}
    \label{fig:vorDNS}
  }
  \subfigure[]
  {
    \includegraphics[scale=0.5]{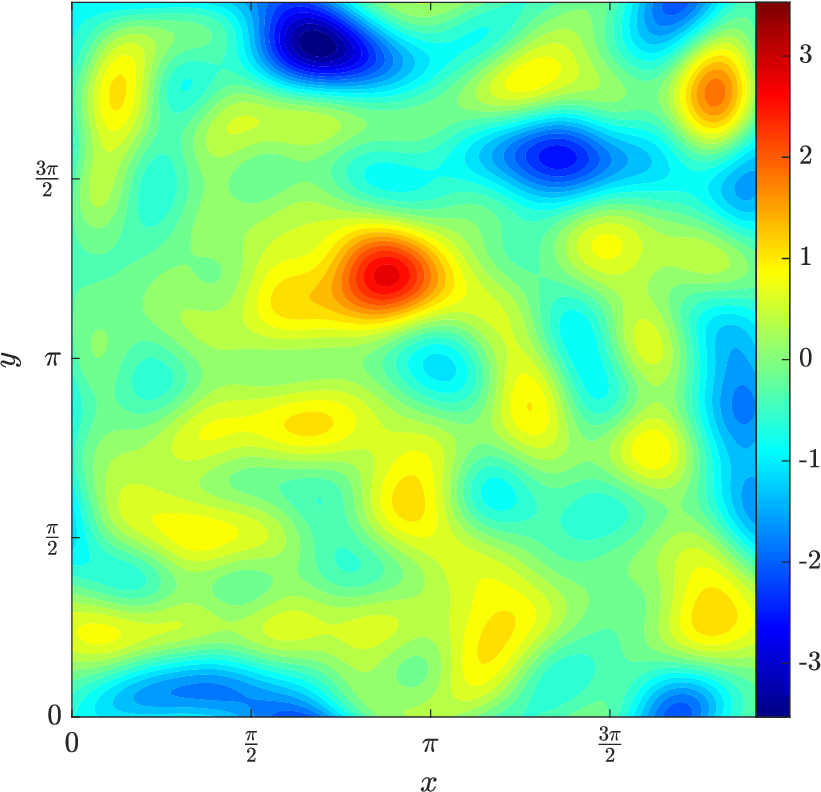}
    \label{fig:vorF}
  }}
\mbox{
  \subfigure[]
  {
    \includegraphics[scale=0.5]{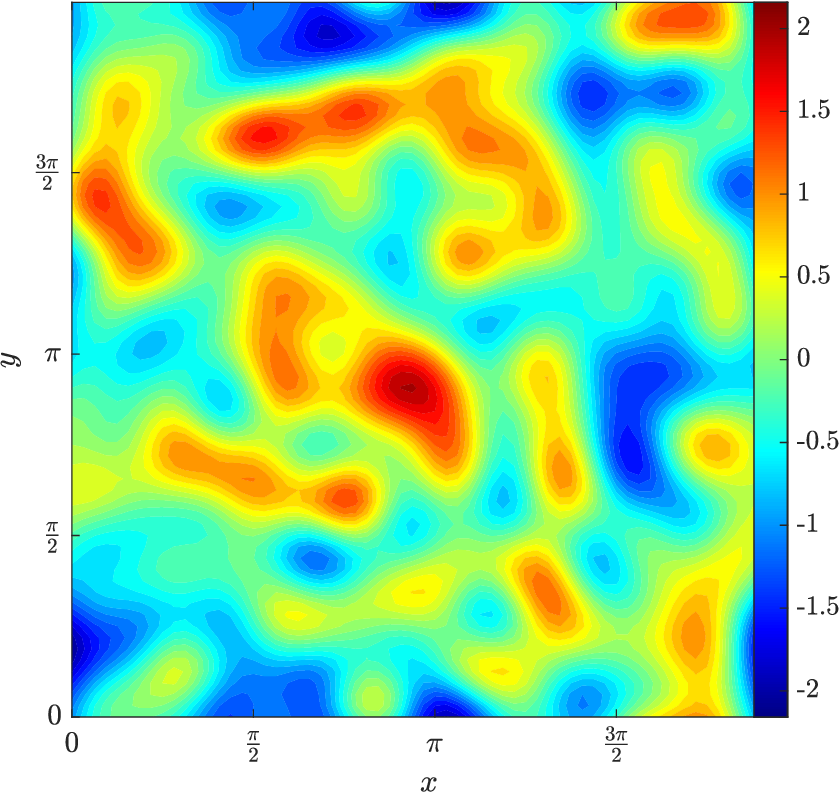}
    \label{fig:vor3A}
  }
  \subfigure[]
  {
    \includegraphics[scale=0.5]{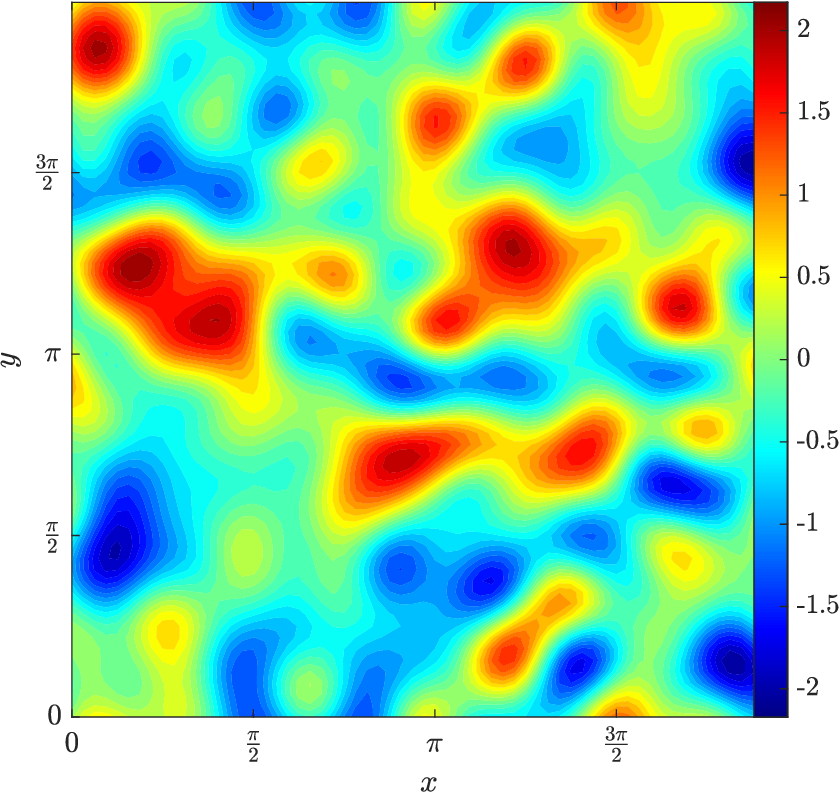}
    \label{fig:vor3B}
  }}
\caption{{(a) The vorticity field $w(T,\x)$ obtained by solving
    the Navier-Stokes system \eqref{eq:2DNS} and (b) the corresponding
    filtered field $\widetilde{w}(T,\x)$. (c,d) The solutions
    $\tw(T,\x)$ of the LES system \eqref{eq:LES} with the optimal eddy
    viscosities $\nuc(s)$ obtained by solving \ref{prob:minJS} and
    \ref{prob:minJM}, respectively. All fields are shown at the time
    instant $t = T$ corresponding to the end of the optimization
    window, cf.~\ref{eq:J2}.}}
\label{fig:vor}
\end{figure}

The optimal eddy viscosities corresponding, via ansatz \eqref{eq:nu},
to the solutions of {Problems \ref{prob:nu}.A and \ref{prob:nu}.B} are
shown as functions of $s \in \I$ in Figure \ref{fig:nuc}. Although
{the eddy viscosities obtained in the two cases {{have
      qualitatively similar profiles and possess common} features such
    as a rapid increase for intermediate values of $\sqrt{s}$}, they
  {also} reveal some fundamentally different physical}
properties. Most notably, the optimal eddy viscosity obtained in
\ref{prob:minJM} is strictly dissipative. On the other hand, this
quantity obtained in the problem in which the constraint \eqref{eq:EE}
is not imposed is negative for small values of the state variable $s$,
meaning that the closure model in fact injects energy in some regions
of the flow. {Such behaviour has {already} been observed
  in closure models \cite{Maulik2020} and was discussed in detail in
  \cite{Matharu2022a}, where optimal closure models {with such
    properties were obtained using different} formulations.} {In
  Figures \ref{fig:vor}a and \ref{fig:vor}b we show, respectively, the
  vorticity field $w(T,\x)$ obtained by solving the Navier-Stokes
  system \eqref{eq:2DNS} and the corresponding filtered field
  $\widetilde{w}(T,\x)$ which serves as the ``target'',
  cf.~\eqref{eq:J2}. As expected, the filtered field reveals fewer
  small-scale features. The solutions $\tw(T,\x)$ of the LES system
  \eqref{eq:LES} with the optimal eddy viscosities $\nuc(s)$ obtained
  by solving \ref{prob:minJS} and \ref{prob:minJM} are shown,
  respectively, in Figures \ref{fig:vor}c and \ref{fig:vor}d. All
  fields are shown at the time instant $t = T$ corresponding to the
  end of the optimization window. Since in the formulation of
  \ref{prob:minJS} and \ref{prob:minJM} matching involves quantities
  defined globally, rather than in a pointwise sense (as was done in
  \cite{Matharu2022a}), the LES fields in Figures \ref{fig:vor}c and
  \ref{fig:vor}d do not appear well correlated with the filtered
  target field in Figure \ref{fig:vor}b. However, they are constructed
  to have similar statistical properties in terms of their enstrophy
  and palinstrophy, cf.~\eqref{eq:E}--\eqref{eq:P}.}

\begin{figure}\centering
  \subfigure[]
  {
    \includegraphics[scale=0.6]{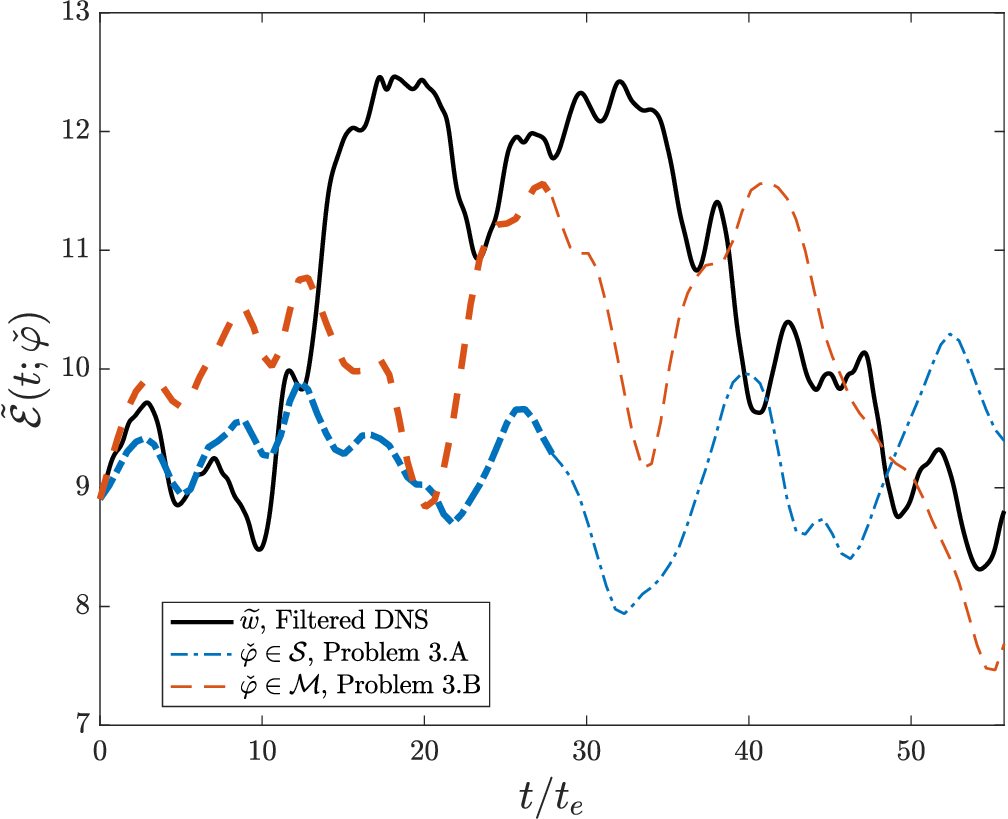}
    \label{figConstraint:Enst}
  }\quad
  \subfigure[]
  {
    \includegraphics[scale=0.6]{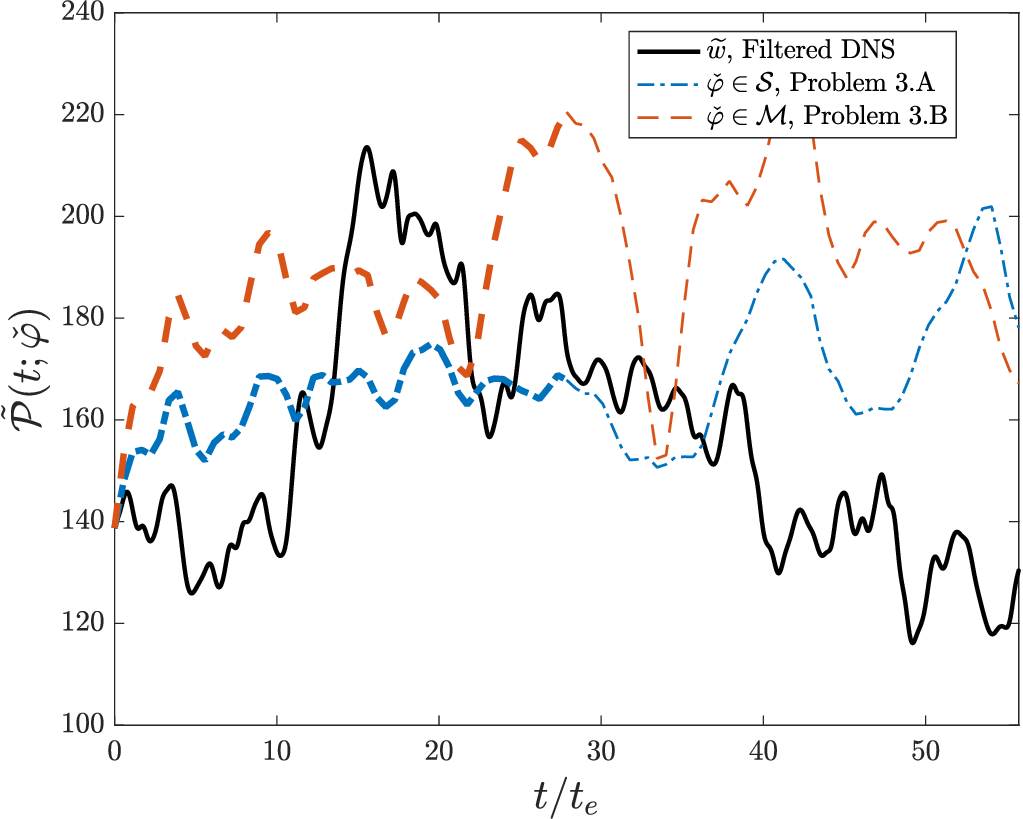}
    \label{figConstraint:Palin}
  }
  \caption{[Problem \ref{prob:nu}:] Dependence of (a) the enstrophy
    \cref{eq:E} and (b) the palinstrophy \cref{eq:P} on time for
    (black, solid line) the filtered DNS flow and the LES flows with
    the optimal eddy viscosities corresponding to (blue, dashed-dot
    line) $\widecheck{{\varphi}} \in \cS$ {found by solving
      Problem \ref{prob:nu}.A} and (red, dashed line)
    $\widecheck{{\varphi}} \in \M$ {found by solving Problem
      \ref{prob:nu}.B}. Thick and thin lines represent to,
    respectively, time in the ``training window'' ($t \in [0, T]$) and
    in the extended window ($t \in (T, 2T]$).}
  \label{fig:tEP}
\end{figure}

To close this section, we show the enstrophy $\tE(t)$ and the
palinstrophy $\tP(t)$ in the filtered DNS flow and in the LES flows
with the optimal eddy viscosities obtained by solving {Problems
  \ref{prob:nu}.A and \ref{prob:nu}.B}, respectively, in Figures
\ref{fig:tEP}a and \ref{fig:tEP}b. These {figures illustrate} how the
matching errors in these two quantities are distributed in time and
are presented on a time window twice longer than the ``training''
window $T$ considered in the error functional \eqref{eq:J2} and the
constraint \eqref{eq:EE}. {We see that for $t \in [0,T]$ the
  behavior of the two quantities is consistent with what was observed
  in Figures \ref{fig:JE3}a and \ref{fig:JE3}b. That is, the better
  reduction of the functional {we see} in Figure~\ref{fig:JE3}a
  results in {an on average} better match of the time history
  of palinstrophy {evident} in Figure~\ref{fig:tEP}b
  {and, conversely, a worse match of the average} enstrophy
  noted in Figure~\ref{fig:tEP}a {is consistent with} deviation
  from the constraint manifold {visible} in
  Figure~\ref{fig:JE3}b. However,} for $t \in [T,2T]$ both enstrophy
and palinstrophy in the LES flows deviate more significantly from
their values in the filtered DNS flow.

\clearpage
\section{Discussion \& Conclusions \label{sec:final}}

{Adjoint-based methods have become a workhorse in the solution of
  unconstrained PDE optimization problems \cite{g03}. They make it
  possible to conveniently determine the gradient (sensitivity) of the
  objective functional with respect to a distributed control variable,
  which can then be used in various gradient descent algorithms
  \cite{nw00}. {Unlike most constraints imposed on the control
    variable, constraints on the state variables are generally harder
    to satisfy since they define, via solutions of the governing
    system,} complicated manifolds in the space of control variables.
  In the present study we demonstrate how the adjoint-based framework
  can be extended to handle such constraints approximately, which is
  done by constructing a projection of the gradient of the objective
  functional onto a subspace tangent to the constraint manifold at the
  given iteration. This projection is realized by solving another
  adjoint problem which is defined in terms of the same adjoint
  operator as the adjoint system employed to determine the gradient,
  but with different {forcing}. Thus, for a state constraint
  defining a codimension-$m$ manifold, a single iteration of the
  gradient algorithm \eqref{eq:desc} requires a total of $(m+1)$
  adjoint solves to evaluate the gradient and the projection on the
  tangent subspace. We focus on the ``optimize-then-discretize''
  paradigm {\cite{g03} in} the infinite-dimensional setting,
  where it is more difficult to use ``black-box'' optimization
  software that could otherwise facilitate imposition of such
  constraints.  Regularity of both the gradient and the projection
  {is} carefully considered, which is done by leveraging the
  Riesz theorem.

  Our approach is illustrated with two test problems, both involving a
  nonlinear scalar constraint: a simple problem describing heat
  conduction in 1D which is used to introduce the method in a clean
  setting and a more complicated problem concerning the optimal design
  of turbulence closures in 2D incompressible flows governed by the
  Navier-Stokes system. This latter problem has a nonstandard
  structure \cite{Matharu2022a} \revtt{and is characterized by a much
    higher level of technical complexity typical of problems arising
    in turbulence research.}  For each of the optimization problems we
  consider the constrained and unconstrained version with and without
  the nonlinear constraint. In addition, for comparison, we also
  consider a version of the first problem where the constraint is
  homogeneous and hence can be enforced exactly with a retraction
  operator.

  Our computation of the gradients and of the normal elements defining
  the tangent subspaces is carefully validated in Section
  \ref{sec:numer}. The computational results presented in Section
  \ref{sec:results} show that even though in the solutions of Problems
  \ref{prob:Heat2}.B and \ref{prob:nu}.B the constraints are not
  satisfied exactly, the obtained optimizers lie much closer to the
  constraint manifold $\M$ than in Problems \ref{prob:Heat2}.A and
  \ref{prob:nu}.A where the constraint is not imposed. However, as
  expected, the optimizers in the constrained problems correspond to
  larger values of the objective functionals. \revtt{In this context,
    we note that it is not a priori known whether Problems
    \ref{prob:Heat2}.B and \ref{prob:nu}.B admit minimizers belonging
    to the constraint manifold for which the objective functional
    would vanish, or nearly vanish (it is, in fact, rather unlikely
    that this could happen). Therefore, a modest reduction of the
    objective functional achieved in these problems, cf.~Figures
    \ref{fig:JE2}a and \ref{fig:JE3}a, should be viewed as a
    reflection of the difficulty of these problems rather than of
    limitations of the proposed approach. In practical applications
    motivating this study even a modest reduction of the objective
    functional is beneficial if the optimal solution satisfies the
    required constraints.}  Importantly, enforcement of the constraint
  via projection onto the subspace tangent to the constraint manifold
  is performed in just as straightforward manner as the computation of
  the gradient.  Therefore, the proposed approach is likely to benefit
  practical applications involving PDE optimization problems.

}


\section*{Acknowledgements}
 
{Partial support for this research was provided through an NSERC (Canada) Discovery Grant: RGPIN-2020-05710. The first author is supported by Vergstiftelsen.}



\bibliographystyle{plainnat}

\end{document}